\documentclass[12pt, reqno]{amsart}
\usepackage{ifthen,amsfonts,amsmath,amssymb,epic,eepic}
\usepackage{epsfig,color}

\makeatother

\theoremstyle{definition}

\def\fnum{equation}
\newtheorem{Thm}[\fnum]{Theorem}
\newtheorem{Cor}[\fnum]{Corollary}

\newtheorem{Lem}[\fnum]{Lemma}

\newtheorem{Def}[\fnum]{Definition}

\newtheorem{Pro}[\fnum]{Proposition}

\newcommand{\Vol}{{\text{Vol}}}

\newcommand{\nn}{{\bf{n}}}

\newcommand{\dist}{{\text {dist}}}

\newcommand{\Hess}{{\text {Hess}}}
\def\ZZ{{\bold Z}}
\def\RR{{\bold R}}

\def\CC{{\bold C }}

\newcommand{\e}{{\text {e}}}

\newcommand{\Genus}{{\text {gen}}}
\newcommand{\cF}{{\mathcal{F}}}

\newcommand{\cD}{{\mathcal{D}}}
\newcommand{\K}{{\text{K}}}
\newcommand{\cL}{{\mathcal{L}}}
\newcommand{\cB}{{\mathcal{B}}}

\newcommand{\cP}{{\mathcal{P}}}
\newcommand{\cS}{{\mathcal{S}}}

\newcommand{\cA}{{\mathcal{A}}}

\newcommand{\rc}{r_1}
\newcommand{\ec}{\epsilon_1}

\newcommand{\eqr}[1]{(\ref{#1})}

\newcommand{\M}{{\bf{M}}}
\newcommand{\cone}{{\bf{C}}}
\newcommand{\ccg}{{\epsilon_a}}

\newcommand{\ccmf}{C_b}
\newcommand{\cdbl}{C_4}
\newcommand{\cmis}{C_6}

\newcommand{\tpi}{\tilde{\Pi}}

\begin{document}

\title[Locally simply connected]
{The space of embedded minimal surfaces of fixed genus in a $3$-manifold IV;
Locally simply connected}

\author{Tobias H. Colding}%
\address{Courant Institute of Mathematical Sciences and Princeton University\\
251 Mercer Street\\ New York, NY 10012 and Fine Hall, Washington Rd.,
Princeton, NJ 08544-1000}
\author{William P. Minicozzi II}%
\address{Department of Mathematics\\
Johns Hopkins University\\
3400 N. Charles St.\\
Baltimore, MD 21218}
\thanks{The first author was partially supported by NSF Grant DMS 9803253
and an Alfred P. Sloan Research Fellowship
and the second author by NSF Grant DMS 9803144
and an Alfred P. Sloan Research Fellowship.}


\email{colding@cims.nyu.edu and minicozz@math.jhu.edu}


\maketitle


\numberwithin{equation}{section}

\section{Introduction} \label{s:s0}

This paper is the fourth in a series where we
describe the space of all embedded minimal surfaces of
fixed genus in a fixed (but arbitrary)
closed $3$-manifold. The key is to understand the structure of an
embedded minimal
disk in a ball in $\RR^3$.
This was undertaken in \cite{CM3}, \cite{CM4} and
the global version of it will be completed here; see \cite{CM15} for
discussion of the local case and \cite{CM13}, \cite{CM14} where we
have surveyed our results about embedded minimal disks.

Our main results are Theorems \ref{t:t0.1}, \ref{t:t2} below.

\begin{figure}[htbp]
    \setlength{\captionindent}{20pt}
    \begin{minipage}[t]{0.5\textwidth}
    \centering\input{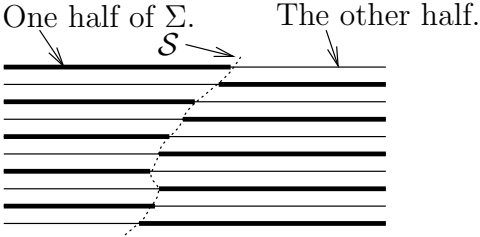}
    \caption{Theorem \ref{t:t0.1} - the singular set, $\cS$, and
the two multi-valued graphs.}
    \end{minipage}
\end{figure}

\begin{Thm} \label{t:t0.1}
See fig. 1.
Let $\Sigma_i \subset B_{R_i}=B_{R_i}(0)\subset \RR^3$ be a
sequence of embedded minimal
disks with $\partial \Sigma_i\subset \partial B_{R_i}$
where $R_i\to \infty$. If $\sup_{B_1\cap \Sigma_i}|A|^2\to \infty$, then
there exists a subsequence, $\Sigma_j$,
and
a Lipschitz curve $\cS:\RR\to \RR^3$ such that after a rotation of $\RR^3$:\\
\underline{1.} $x_3(\cS(t))=t$.  (That is, $\cS$ is a graph
over the $x_3$-axis.)\\
\underline{2.}  Each $\Sigma_j$ consists of exactly two multi-valued
graphs away
from $\cS$ (which spiral together).\\
\underline{3.} For each $1>\alpha>0$, $\Sigma_j\setminus \cS$ converges
in the $C^{\alpha}$-topology to the foliation,
$\cF=\{x_3=t\}_t$, of $\RR^3$.\\
\underline{4.}  For all $r>0$, $t$, then
$\sup_{B_{r}(\cS (t))\cap \Sigma_j}|A|^2\to\infty$.
\end{Thm}

In \underline{2.}, \underline{3.}
that $\Sigma_j\setminus \cS$ are multi-valued graphs and
converges to $\cF$ means that for each compact subset
$K\subset \RR^3\setminus \cS$
and $j$ sufficiently large $K\cap \Sigma_j$ consists of multi-valued
graphs over (part of) $\{x_3=0\}$
and $K\cap \Sigma_j\to K\cap \cF$.

This theorem (as many of the results below)
is modeled by the helicoid and its rescalings.
The helicoid is
the minimal surface $\Sigma^2$ in $\RR^3$ given by
$(s\,\cos t,s\sin t,-t)$
where $s,t\in\RR$. Take a sequence
$\Sigma_i = a_i \, \Sigma$ of rescaled helicoids where
$a_i \to 0$. Since the helicoid has cubic volume growth,
the density is unbounded. The curvature is blowing up
along the vertical axis. The sequence converges
(away from the vertical axis)
to a foliation by flat parallel planes.
The singular
set $\cS$ (the axis) then consists of removable singularities.

\begin{figure}[htbp]
    \setlength{\captionindent}{20pt}
\begin{minipage}[t]{0.5\textwidth}
    \centering\input{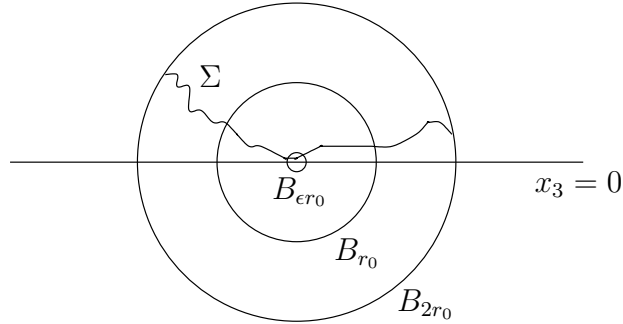}
    \caption{Theorem \ref{t:t2} - the one-sided curvature estimate for an
embedded minimal disk $\Sigma$ in a half-space with 
$\partial \Sigma\subset \partial B_{2r_0}$:  The components of
$B_{r_0}\cap \Sigma$ intersecting $B_{\epsilon r_0}$ are graphs.}
    \end{minipage}
\end{figure}

\begin{Thm}  \label{t:t2}
See fig. 2.  There exists $\epsilon>0$, such that if
$\Sigma^2\subset B_{2r_0} \cap \{x_3>0\}
\subset \RR^3$ is an
embedded minimal
disk with $\partial \Sigma\subset \partial B_{2 \, r_0}$,
then for all components $\Sigma'$ of
$B_{r_0} \cap \Sigma$ which intersect $B_{\epsilon r_0}$
\begin{equation}	\label{e:graph}
\sup_{\Sigma'} |A_{\Sigma}|^2
\leq r_0^{-2} \, .
\end{equation}
\end{Thm}

Using the minimal surface equation and that $\Sigma'$ has points
close to a plane, it is not hard to see that, for $\epsilon>0$
sufficiently small, \eqr{e:graph} is equivalent to the statement 
that $\Sigma'$
is a graph over the plane $\{x_3=0\}$.

An embedded minimal surface $\Sigma$ which is as in
Theorem \ref{t:t2} is said to satisfy the
$(\epsilon , r_0)$-{\it effective one-sided Reifenberg condition};
cf. Appendix \ref{s:ref}.  We will often refer to Theorem \ref{t:t2}
as {\it the one-sided curvature estimate}.

\begin{figure}[htbp]
    \setlength{\captionindent}{20pt}
    \begin{minipage}[t]{0.5\textwidth}
    \centering\input{unot6.pstex_t}
    \caption{The catenoid given by revolving $x_1= \cosh x_3$
around the $x_3$-axis.}
    \end{minipage}\begin{minipage}[t]{0.5\textwidth}
    \centering\input{unot7.pstex_t}
    \caption{Rescaling the catenoid shows that simply connected is
    needed in the one-sided curvature estimate.}
    \end{minipage}%

\end{figure}

Note that the assumption in Theorem \ref{t:t2}
that $\Sigma$ is simply connected is crucial
as can be seen from the example of a rescaled catenoid. The catenoid
is the minimal surface in $\RR^3$ given by
$(\cosh s\, \cos t,\cosh s\, \sin t,s)$
where $s,t\in\RR$; see fig. 3.
Under rescalings this converges (with multiplicity two) to
the flat plane; see fig. 4.

An almost immediate consequence of Theorem \ref{t:t2} is:

\begin{figure}[htbp]
    \setlength{\captionindent}{20pt}
    \begin{minipage}[t]{0.5\textwidth}
    \centering\input{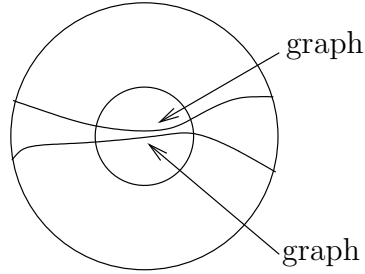}
    \caption{Corollary \ref{c:barrier}:  Two sufficiently close
    components of an embedded minimal disk must each be a graph.}
    \end{minipage}
\end{figure}

\begin{Cor}      \label{c:barrier}
See fig. 5.  There exist $c > 1$, $\epsilon >0$ so:
Let $\Sigma_1,\, \Sigma_2 \subset B_{cr_0} \subset \RR^3$ be
disjoint embedded minimal
surfaces with $\partial \Sigma_i \subset \partial B_{cr_0}$
and $B_{\epsilon \, r_0 } \cap \Sigma_i \ne \emptyset$.
If $\Sigma_1 $ is a disk,
 then for all components $\Sigma_1'$ of
$B_{r_0} \cap \Sigma_1$ which intersect
$B_{\epsilon \, r_0}$
\begin{equation}        \label{e:onece}
    \sup_{\Sigma_1'}   |A|^2
        \leq  r_0^{-2}  \, .
\end{equation}
\end{Cor}

To explain how these theorems are proven using the results of
\cite{CM3}--\cite{CM5}, and \cite{CM7} we will need some notation
for multi-valued graphs.  Let $\cP$ be the universal cover of the
punctured plane $\CC \setminus \{ 0 \}$ with global (polar)
coordinates $(\rho , \theta)$ and set $S_{r,s}^{\theta_1 ,
\theta_2} = \{  r \leq  \rho \leq  s  \, , \, \theta_1 \leq \theta
\leq \theta_2 \}$. An $N$-valued graph $\Sigma$ of a function $u$
over the annulus $D_{s} \setminus D_{r}$ is a
(single-valued) graph (of $u$) over $S_{r,s}^{-N\,\pi,N\,\pi}$
($\Sigma_{r,s}^{\theta_1 , \theta_2}$ will denote the subgraph of
$\Sigma$ over $S_{r,s}^{\theta_1 , \theta_2}$). The separation, see fig. 6,
between consecutive sheets will be denoted by $w$ so
$w(\rho,\theta) = u(\rho,\theta + 2 \, \pi) - u(\rho,\theta)$.  A
multi-valued graph is embedded if and only if $|w|>0$. Note (see fig. 7) that
one-half of the helicoid, i.e.,  each of the two components of
$(s\cos t,s\sin t,-t)\setminus \{s=0\}$, as an $\infty$-valued
graph of a function given in polar coordinates by
$u(\rho,\theta)=-\theta$, $u(\rho,\theta)=-\theta+\pi$,
respectively.  In particular, $w(\rho,\theta)=-2\,\pi$.

In this paper, as in \cite{CM7}, we have normalized so embedded 
multi-valued graphs have negative seperation.  This can be achieved 
after possibly reflecting in a plane.  

\begin{figure}[htbp]
    \setlength{\captionindent}{20pt}
    \begin{minipage}[t]{0.5\textwidth}
    \centering\input{unot2a.pstex_t}
    \caption{The separation of a multi-valued graph.
	(Here the multi-valued graph is shown with
	negative separation.)}
    \end{minipage}\begin{minipage}[t]{0.5\textwidth}
    \centering\input{unot1.pstex_t}
    \caption{The helicoid is obtained
    by gluing together two $\infty$-valued graphs along a line.  
The two multi-valued graphs are given in polar coordinates
by $u_1(\rho,\theta)=-\theta$ and $u_2(\rho,\theta)=-\theta+\pi$.  
In either case $w(\rho,\theta)=-2\,\pi$.}
    \end{minipage}%
\end{figure}

In \cite{CM4} we showed that an embedded minimal disk in a ball
in $\RR^3$ with large curvature at a point contains an
almost flat multi-valued graph nearby.  Namely, we showed:

\begin{figure}[htbp]
    \setlength{\captionindent}{20pt}
    \begin{minipage}[t]{0.5\textwidth}
    \centering\input{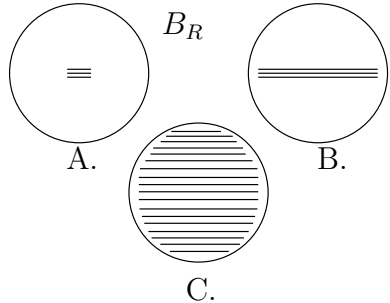}
    \caption{Proving Theorem \ref{t:t0.1}.  A. Finding a small
$N$-valued graph in $\Sigma$.
    B. Extending it in $\Sigma$ to a large
    $N$-valued graph.  C. Extend
    the number of sheets.  (A. follows from \cite{CM4}
    and B. follows from \cite{CM3}.)}
    \end{minipage}
\end{figure}

\begin{Thm} \label{t:blowupwinding0}
(Theorem $0.2$ of \cite{CM4}).  See A. and B. in fig. 8.  
Given $N\in \ZZ_+$, $\epsilon > 0$, there exist
$C_1,\,C_2>0$ so: Let
$0\in \Sigma^2\subset B_{R}\subset \RR^3$ be an embedded minimal
disk, $\partial \Sigma\subset \partial B_{R}$. If
$\sup_{B_{r_0}}|A|^2\geq 4\,C_1^2\,r_0^{-2}$ for some
$0<r_0<R$, then there exists
(after a rotation)
an $N$-valued graph $\Sigma_g \subset \Sigma\cap
\{ x_3^2 \leq \epsilon^2 \, (x_1^2 + x_2^2) \}$ over $D_{R/C_2}
\setminus D_{2r_0}$ with gradient $\leq \epsilon$.
\end{Thm}

An important consequence of Theorem \ref{t:blowupwinding0} is (see
theorem 5.8 of \cite{CM4}): Let $\Sigma_i\subset B_{2R}$ be a
sequence of embedded minimal
 disks with
$\partial \Sigma_i\subset \partial B_{2R}$.  Clearly
(after possibly going to a subsequence) either (1) or (2) occur:\\
(1) $\sup_{B_{R}\cap\Sigma_i}|A|^2\leq C<\infty$ for some constant $C$.\\
(2) $\sup_{B_{R}\cap\Sigma_i}|A|^2\to \infty$.\\
In (1) (by a standard argument) $\cB_s(y_i)$ is a
graph for all $y_i\in B_{R}\cap \Sigma_i$, where $s$ depends only on $C$.
In (2) (by theorem 5.8 of \cite{CM4}) if
$y_i\in B_{R}\cap \Sigma_i$ with $|A|^2(y_i)\to \infty$, then we can after
passing to a subsequence assume that $y_i\to y$, each $\Sigma_i$ contains a
$2$-valued graph $\Sigma_{d,i}$ over
$D_{R/C_2}(y)\setminus D_{\epsilon_i}(y)$ with
$\epsilon_i\to 0$, and
$\Sigma_{d,i}$ converges to a graph $y\in \Sigma_{\infty}$ over $D_{R/C_2}(y)$.
In either
case in the limit there is a smooth minimal graph through each point in the
support.

These multi-valued graphs should be thought of as the basic
building blocks for an embedded minimal disk. In fact, using a
standard blow up argument, we showed in \cite{CM4} (corollary
$4.14$ combined with proposition $4.15$ there) that Theorem
\ref{t:blowupwinding0} was a consequence of the following that we
will use to construct the actual building blocks starting off on
the smallest possible scale:

\begin{Thm} \label{t:cm34}
\cite{CM4}. Given $N\in \ZZ_+$, $\epsilon > 0$, there exist
$C_1,\,C_2, C_3>0$ so: Let $0\in \Sigma^2\subset B_{R}\subset
\RR^3$ be an embedded minimal disk, $\partial \Sigma\subset
\partial B_{R}$. If $\sup_{B_{r_0} \cap \Sigma}|A|^2\leq
4\,C_1^2\,r_0^{-2}$ and $|A|^2(0)=C_1^2\,r_0^{-2}$ for some
$0<r_0<R$, then there exists (after a rotation) an $N$-valued
graph $\Sigma_g \subset \Sigma\cap \{x_3^2\leq \epsilon^2 \,
(x_1^2 + x_2^2) \}$ over $D_{R/C_2} \setminus D_{r_0}$ with
gradient $\leq \epsilon$ and separation $\geq C_3\,r_0$ over
$\partial D_{r_0}$.
\end{Thm}

It will be important for the application of Theorem \ref{t:cm34}
here that the initial separation of the sheets
is proportional to the initial scale that the graph starts off on.

Theorems \ref{t:t0.1} and \ref{t:t2} deal with how the building blocks
fit together.
A consequence of Theorem \ref{t:t0.1} is
that if an embedded minimal disk starts to spiral very tightly, then it
can unwind only
very slowly.  That is, in a whole extrinsic tubular neighborhood it
continues to spiral tightly and fills up almost the entire space.

Let us briefly outline the proof of the one-sided; i.e.,
Theorem \ref{t:t2}. Suppose that $\Sigma$
is an embedded minimal disk in the half-space $\{x_3> 0\}$. We prove the
curvature estimate by contradiction; so suppose that $\Sigma$ has low
points with large curvature.
Starting at such a point, we decompose
(see Corollary \ref{c:decompo})
$\Sigma$ into disjoint multi-valued graphs using the existence
of nearby points with large curvature (see Proposition \ref{c:0.4ofcm5}),
 a blow up argument,
and \cite{CM3}, \cite{CM4}. The key point is then to show
(see Proposition \ref{p:lift} and fig. 9)
that we can in fact find such a decomposition where the ``next''
multi-valued graph starts off a definite amount below where the
previous multi-valued graph started off. In fact, what we show
is that this definite amount is a fixed fraction of the distance
between where the two graphs started off. Iterating this
eventually forces $\Sigma$ to have points where $x_3<0$.
Which is the desired contradiction.

\begin{figure}[htbp]
    \setlength{\captionindent}{20pt}
    \begin{minipage}[t]{0.5\textwidth}
    \centering\input{unot10.pstex_t}
    \caption{Two consecutive blow up points satisfying \eqr{e:pairs}.}
    \end{minipage}
\end{figure}

\begin{figure}[htbp]
    \setlength{\captionindent}{20pt}
    \begin{minipage}[t]{0.5\textwidth}
    \centering\input{unot11.pstex_t}
    \caption{Between two consecutive blow up points satisfying  \eqr{e:pairs}
there are a bunch of blow up points satisfying
 \eqr{e:defc1ii}.}
    \end{minipage}%
    \begin{minipage}[t]{0.5\textwidth}
    \centering\input{unot12.pstex_t}
    \caption{Measuring height.  Blow up points and
corresponding multi-valued graphs.}
    \end{minipage}
\end{figure}

To prove this key proposition
(Proposition \ref{p:lift}) we use two decompositions
and two kinds of blow up points.
The first decomposition which is Corollary \ref{c:decompo} uses
the more standard blow up points given by \eqr{e:defc1ii}.  These
are  pairs $(y,s)$
where $y\in \Sigma$ and $s>0$ is such that
$\sup_{\cB_{8s}(y)}|A|^2\leq 4|A|^2(y)=4C_1^2s^{-2}$.
The point about such a pair $(y,s)$ is that by \cite{CM3}, \cite{CM4}
(and an argument in
Part \ref{s:sotherhalf}
which allows us replace extrinsic balls by
intrinsic ones), then
$\Sigma$ contains a multi-valued graph near $y$ starting
off on the scale $s$.  (This is assuming that $C_1$ is a sufficiently
large constant given by \cite{CM3}, \cite{CM4}.)  The second kind of blow up
points are the ones satisfying \eqr{e:pairs}.  Basically \eqr{e:pairs}
is \eqr{e:defc1ii}
(except for a technical issue)  where $8$ is replaced by some really large
constant $C$.  The point will then be that we can find blow up points
satisfying \eqr{e:pairs} so that the distance between them is proportional
to the sum of the scales.  Moreover, between consecutive
blow up points satisfying \eqr{e:pairs}, we
can find a bunch of blow up points satisfying \eqr{e:defc1ii}; 
see fig. 10.  The advantage
is now that if we look between blow up points satisfying \eqr{e:pairs}, then
the height of the multi-valued graph given by such a pair grows
like a small power of the distance whereas the separation between the sheets
in a multi-valued graph given by \eqr{e:defc1ii} decays like a small
power of the
distance; see fig. 11.  
Now thanks to that the number of blow up points satisfying
\eqr{e:defc1ii} (between two consecutive blow up points satisfying
\eqr{e:pairs}) grows
almost linearly then, even though the height of the graph coming from
the blow up point satisfying \eqr{e:pairs} could move up (and thus work
against us), then the sum of the separations of the graphs coming from the
points satisfying \eqr{e:defc1ii} dominates the other term.  Thus the
next blow
up point satisfying \eqr{e:pairs} (which lies below all the other graphs)
is forced to be a definite amount lower than the previous blow up
point satisfying \eqr{e:pairs}.

Let $x_1 , x_2 , x_3$ be the standard coordinates on $\RR^3$ and
$\Pi : \RR^3 \to \RR^2$ orthogonal projection to $\{ x_3 = 0 \}$.
For $y \in S \subset \Sigma \subset \RR^3$ and $s > 0$, the
extrinsic and intrinsic balls are
$B_s(y)$, $\cB_s(y)$ and
$\Sigma_{y,s}$ is the component of $B_{s}(y)
\cap \Sigma$ containing $y$.
$D_s$ denotes the disk $B_s(0) \cap \{ x_3 = 0 \}$.
$\K_{\Sigma}$ the sectional curvature of a smooth compact surface
$\Sigma$ and when
$\Sigma$ is immersed $A_{\Sigma}$ will be its second fundamental form.
When $\Sigma$ is oriented, $\nn_{\Sigma}$ is the unit normal.

This paper completes the results announced in \cite{CM11} and
\cite{CM12}.

Using Theorems \ref{t:t0.1}, \ref{t:t2},
W. Meeks and H. Rosenberg proved
that the plane and helicoid are the only
complete properly embedded
simply-connected minimal surfaces in $\RR^3$, \cite{MeRo}.


\setcounter{part}{0}
\numberwithin{section}{part} 
\renewcommand{\rm}{\normalshape} 
\renewcommand{\thepart}{\Roman{part}}
\setcounter{section}{1}

\part{The proof of Theorem \ref{t:t0.1} assuming Theorem
\ref{t:t2} and short curves} \label{p:p2}

In this part we will show how Theorem \ref{t:t0.1} follows from
Theorem \ref{t:t2}, the results about existence of multi-valued
graphs from \cite{CM3}, \cite{CM4} which were recalled in the
introduction, corollary III.3.5 of \cite{CM5}, and the results
about properness of embedded disks from \cite{CM7} (once we see
that the conditions in corollary $0.7$ of \cite{CM7} are
satisfied). The remaining parts of this paper are devoted to
showing Theorem \ref{t:t2} (Part \ref{p:p00})
and that corollary $0.7$ of \cite{CM7} applies (Part
\ref{p:shc}; see, in particular, Theorem \ref{t:main}
below).

We will use several times that given $\alpha > 0$, Proposition
II.2.12 of \cite{CM3} gives $N_g$ so if $u$ satisfies the minimal
surface equation on $S_{\e^{-N_g}, \e^{N_g} \, R}^{- N_g, 2 \pi +
N_g}$ with $|\nabla u| \leq 1$, and $w<0$, then $\rho \, |\Hess_u|
+ \rho \, |\nabla  w |/|w| \leq \alpha$ on $S_{1, R}^{0 , 2\pi }$.
Theorem 3.36 of \cite{CM9} then yields $|\nabla u - \nabla u(1,0)|
\leq C \alpha$.  We can therefore assume (after rotating so
$\nabla u(1,0) = 0$) that
\begin{equation} \label{e:wantit}
    |\nabla u| + \rho \, |\Hess_u| +  4 \, \rho \, |\nabla  w |/|w| +
    \rho^2 \, |\Hess_w  |/|w|
      \leq \epsilon < 1/(2\pi) \, .
\end{equation}
The bound on $|\Hess_w  |$ follows from the other bounds and
standard elliptic theory.
 In what follows, we will assume that $w
< 0$. (This normalizes the graph of $u$ to spiral downward; this
can be achieved after possibly reflecting in a plane.)

If $\Sigma$ is an embedded graph of $u$ over
$S_{1/2,2R}^{-3\pi,N+3\pi}$, then $E$ is the region over $D_{R}
\setminus D_{1}$ between the top and bottom sheets of the
concentric  subgraph over $S_{1,R}^{-2\pi,N+2\pi}$ (recall that,
possibly after reflection, we can assume $w< 0$).  Namely, when
$N$ is even, $E$ is the set (see fig. 12)
of all $(r \cos \theta , r \sin \theta
, t)$ with $1 \leq r \leq R$, $-2 \pi \leq \theta < 0$, and
\begin{equation}    \label{e:defEi}
    u (r , \theta + (N+2) \pi) < t < u (r , \theta )
     \, .
\end{equation}

To apply corollary $0.7$ of \cite{CM7} we need the following
result (which will be proven in Part \ref{p:shc}) on
existence of ``the other half'' of an embedded minimal disk and
short curves, $\sigma_{\theta}$, connecting the two halves:

\begin{figure}[htbp]
    \setlength{\captionindent}{20pt}
    \begin{minipage}[t]{0.5\textwidth}
    \centering\input{uls12.pstex_t}
    \caption{The set $E$ in \eqr{e:defEi}.}
    \end{minipage}\begin{minipage}[t]{0.5\textwidth}
    \centering\input{uls13.pstex_t}
    \caption{Theorem \ref{t:main} - the existence of the 
other half and the short curves, $\sigma_{\theta}$, connecting the two halves.}
    \end{minipage}
\end{figure}

\begin{Thm}  \label{t:main}
See fig. 13.
There exist $C$, $R_0$, $N_0$, $\epsilon > 0$ so for $N\geq N_0$: Let
$\Sigma \subset B_{4R}$ be an embedded minimal disk, $\partial
\Sigma\subset
\partial B_{4R}$, $R\geq R_0$, and
 $\Sigma_{1} \subset \Sigma$ a graph
of $u_1$ with $|\nabla u_1| \leq \epsilon$ over
$S_{1/2,2R}^{-3\pi,N+3\pi}$.
  Then $E  \cap \Sigma \setminus \Sigma_{1}$ is
a graph  of $u_2$ over $S_{1,R}^{0,N+2\pi}$ and $u_1(1,2\pi) < u_2
(1,0) < u_1 (1,0)$. Moreover, for all $0 \leq \theta \leq N
+2\pi$, a curve $\sigma_{\theta}\subset \{ x_1^2+x_2^2\leq 1 \}
\cap \Sigma $ with length $\leq C$ connects the image of $u_1$
over $(1,\theta)$ with the image of $u_2$ over $(1,\theta)$.
\end{Thm}

The main example of the ``two halves'' of an embedded minimal disk
and short curves connecting them comes from the helicoid.  Namely,
let $\Sigma$ be the helicoid, i.e., $\Sigma=(\rho\,\cos
\theta,\rho\,\sin \theta, -\theta)$ where $\rho,\,\theta\in\RR$,
then $\Sigma\setminus \{\rho=0\}$ consists of two $\infty$-valued
graphs $\Sigma_1$, $\Sigma_2$ and  $\sigma_\theta$ given by
$\Sigma\cap \{x_3=- \theta\}$ union  $\{ (-\cos \tau,- \sin \tau,
-\tau  ) \,|\,\theta \leq \tau\leq \theta+  \pi\} $ are
short curves connecting the two halves.  Theorem \ref{t:main}
asserts that this is the general picture.

\vskip2mm
We will use the following
result from \cite{CM5} to get nearby points with large curvature
(here, as before, $\Sigma_{y,s}$ is the component of $B_{s}(y)
\cap \Sigma$ containing $y$):

\begin{figure}[htbp]
    \setlength{\captionindent}{20pt}
    \begin{minipage}[t]{0.5\textwidth}
    \centering\input{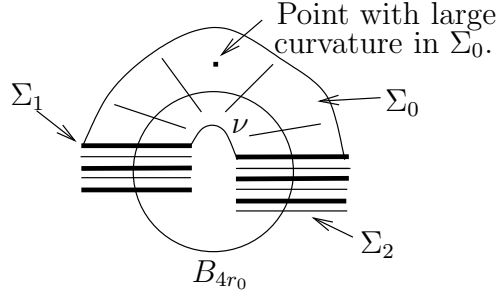}
    \caption{Proposition \ref{c:0.4ofcm5} - 
existence of nearby points with large curvature.}
    \end{minipage}
\end{figure}

\begin{Pro} \label{c:0.4ofcm5}
(Corollary III.3.5 of \cite{CM5}).  See fig. 14.  Given $C_1$, there exists
$C_2$ so: Let  $0 \in\Sigma \subset B_{2 C_2 \, r_0}$ be an
embedded minimal disk. Suppose $ \Sigma_1 , \Sigma_2 \subset
\Sigma\cap \{x_3^2\leq (x_1^2+x_2^2)\}$ are graphs of $u_i$
satisfying \eqr{e:wantit} on $S_{r_0, C_2 r_0}^{-2\pi,2\pi}$,
$u_1(r_0,2\pi) < u_2 (r_0,0) < u_1 (r_0,0)$, and $\nu \subset
\partial \Sigma_{0,2r_0}$ a curve from $\Sigma_1$ to $\Sigma_2$.
Let $\Sigma_0$ be the component of $\Sigma_{0,C_2 r_0} \setminus
(\Sigma_1 \cup \Sigma_2 \cup \nu)$ which does not contain
$\Sigma_{0,r_0}$.  Suppose  either
 $\partial \Sigma \subset \partial B_{2 C_2 \, r_0}$ or
$\Sigma$ is stable and $\Sigma_0$  does not intersect $\partial
\Sigma$. Then
\begin{equation}    \label{e:nextone1c}
    \sup_{x\in \Sigma_{0} \setminus B_{4r_0}} |x|^2 \, |A|^2 (x)
\geq 4\,C_1^2 \, .
\end{equation}
\end{Pro}

Note that by the curvature estimate for stable surfaces, \cite{Sc},
\cite{CM2}, when $\Sigma$ is stable then the conclusion of Proposition
\ref{c:0.4ofcm5} is that no such surface exists for $C_1$,
$C_2$ sufficiently large.

\section{Regularity of the singular set}   \label{s:ref1}

If $\delta>0$ and $z\in \RR^3$, then we denote by
$\cone_{\delta}(z)$ the (convex) cone with vertex $z$, cone angle $(\pi/2 -
\arctan \delta)$, and axis parallel to the $x_3$-axis.    That is, see fig. 15,
\begin{equation}
\cone_{\delta}(z)=\{x\in \RR^3\,|\,(x_3-z_3)^2 \geq
\delta^2\,((x_1-z_1)^2+(x_2-z_2)^2) \} \, .
\end{equation}

\begin{figure}[htbp]
    \setlength{\captionindent}{20pt}
\begin{minipage}[t]{0.5\textwidth}
    \centering\input{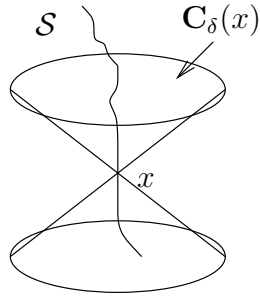}
    \caption{It follows from the one-sided curvature estimate that the singular
    set has the cone property and hence is a Lipschitz curve; see 
    Lemma \ref{l:regsing}.}
    \end{minipage}
\end{figure}

\begin{Lem}  \label{l:regsing}
See fig. 15.  Let $0\in\cS\subset \RR^3$ be a closed set such that for some
$\delta>0$ and each $z\in \cS$, then $\cS\subset
\cone_{\delta}(z)$. If for all $t\in x_3(\cS)$ and all
$\epsilon>0$, $\cS\cap \{t<x_3<t+\epsilon\}\ne \emptyset$,
$\cS\cap \{t-\epsilon<x_3<t\}\ne \emptyset$, then $\cS\cap
\{x_3=t\}$ consists of exactly one point $\cS_t$ for all
$t\in\RR$, and $t\to \cS_t$ is a Lipschitz parameterization of
$\cS$.  In fact,
\begin{equation}  \label{e:lipscbd}
|t_2-t_1|\leq |\cS_{t_2}-\cS_{t_1}|\leq \sqrt{1+\delta^{-2}}\,|t_2-t_1|\, .
\end{equation}
\end{Lem}

\begin{proof}
First by the cone property it follows that $\cS\cap \{x_3=t\}$
consists of at most one point for each $t\in \RR$.  Assume that
$\cS\cap \{x_3=t_0\}=\emptyset$ for some $t_0$.  Since $\cS\subset
\RR^3$ is a nonempty closed set and $x_3:\cS\subset
\cone_{\delta}(0)\to \RR$ is proper, then $x_3(\cS)\subset \RR$ is
also closed (and nonempty). Let $t_s\in x_3(\cS)$ be the closest
point in $x_3(\cS)$ to $t_0$. The desired contradiction now easily
follows since either $\cS\cap \{t_s<x_3<t_0\}$ or $\cS\cap
\{t_0<x_3<t_s\}$ is nonempty by assumption.

It follows that $t\to \cS_t$ is a well-defined curve (from $\RR$
to $\cS$). Moreover, since $\cS_{t_2}\subset
\{x_3=t_1+(t_2-t_1)\}\cap \cone_{\delta}(\cS_{t_1}) \subset
B_{\sqrt{1+\delta^{-2}}|t_2-t_1|}(\cS_{t_1})$, \eqr{e:lipscbd}
follows.
\end{proof}

We will refer loosely to a set $\cS$ as in Lemma \ref{l:regsing}
as having the cone property.
Next we will see, by a very general compactness argument, that for any
sequence of surfaces in $\RR^3$, after possibly going to a subsequence,
then there is a well defined notion of points where the second
fundamental form of the sequence blows up.  The set of such points will
below be the $\cS$ in Lemma \ref{l:regsing}; we observe
in Corollary \ref{c:conecor}  below that $\cS$ has the cone property.

\begin{Lem} \label{l:inftyornot}
Let $\Sigma_i\subset B_{R_i}$, $\partial \Sigma_i\subset \partial B_{R_i}$,
and $R_i\to \infty$ be a sequence of (smooth) compact surfaces.
After passing to a subsequence, $\Sigma_j$, we
may assume that for each $x\in \RR^3$ either $\sup_{B_{r}(x)\cap
\Sigma_j}|A|^2\to \infty$ for all $r>0$ or $\sup_j\sup_{B_r(x)\cap
\Sigma_j}|A|^2<\infty$ for some $r>0$.
\end{Lem}

\begin{proof}
For $r>0$ and an integer $n$, define a sequence of functions on $\RR^3$ by
\begin{equation}
\cA_{i,r,n}(x)=\min \{ n,\sup_{B_r(x)\cap \Sigma_i}|A|^2\}\, ,
\end{equation}
where we set $\sup_{B_r(x)\cap \Sigma_i}|A|^2=0$ if
$B_r(x)\cap \Sigma_i=\emptyset$.
Set
\begin{equation}
\cD_{i,r,n}=\lim_{k\to \infty} 2^{-k} \,
\sum_{m=0}^{2^k-1} \cA_{i,(1+m2^{-k})r,n}\, ,
\end{equation}
then $\cD_{i,r,n}$ is continuous and $\cA_{i,2r,n}\geq \cD_{i,r,n}\geq
\cA_{i,r,n}$.
 Let $\nu_{i,r,n}$
be the (bounded) functionals,
\begin{equation}
\nu_{i,r,n}(\phi)=\int_{B_n}\phi\,\cD_{i,r,n}\text{ for }
\phi\in L^2 (\RR^3)\, .
\end{equation}
By standard compactness for fixed $r$, $n$, after passing to a subsequence,
$\nu_{j,r,n}\to \nu_{r,n}\text{ weakly}$.
In fact (since the unit ball in $L^2(\RR^3)$ has a countable basis),
by an easy diagonal argument after passing to a
subsequence we may assume that for all $n, m \geq 1$ fixed
$\nu_{j,2^{-m},n}\to \nu_{2^{-m},n}\text{ weakly}$.
Note that if $x\in \RR^3$ and for all $m$, $n$ with $n\geq |x|+1$,
(identify $B_{2^{-m}}(x)$ with its characteristic function)
\begin{equation} \label{e:inftyornot}
\nu_{2^{-m},n}(B_{2^{-m}}(x))\geq n\, \Vol (B_{2^{-m}})\, ,
\end{equation}
then for each fixed $r>0$, $\sup_{B_r(x)\cap\Sigma_j}|A|^2\to \infty$.
Conversely, if for some $n \geq |x| +1$, $m$, \eqr{e:inftyornot}
fails at $x$, then
$\sup_j\sup_{B_r(x)\cap
\Sigma_j}|A|^2<\infty$ for $r=2^{-m-1}$.
\end{proof}

To implement Lemma \ref{l:regsing} in the proof of Theorem \ref{t:t0.1},
we will need the following (direct) consequence of Theorem \ref{t:t2} with
$\Sigma_d$ playing the role of the plane
(and the maximum principle as in Appendix \ref{s:p4}):

\begin{figure}[htbp]
    \setlength{\captionindent}{20pt}
    \begin{minipage}[t]{0.5\textwidth}
    \centering\input{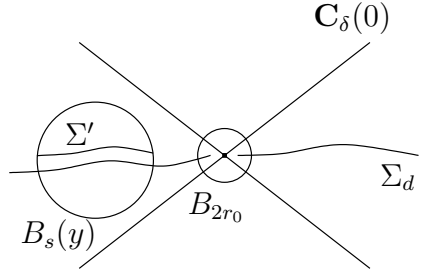}
    \caption{Corollary \ref{c:conecor}:  With $\Sigma_d$ playing
    the role of $x_3=0$,
    by the one-sided estimate, $\Sigma$ consists of multi-valued
    graphs away from a cone.}
    \end{minipage}

\end{figure}

\begin{Cor}   \label{c:conecor}
See fig. 16.  There exists $\delta_0>0$ so: Suppose $\Sigma\subset B_{2R}$,
$\partial \Sigma\subset \partial B_{2R}$ is an embedded minimal
disk containing a $2$-valued graph $\Sigma_d \subset \{x_3^2 \leq
\delta_0^2\, (x_1^2+x_2^2)\}$ over $D_{R}\setminus D_{r_0}$ with
gradient $\leq \delta_0$.  Then each component of $B_{R/2}\cap
\Sigma\setminus (\cone_{\delta_0}(0)\cup B_{2 r_0})$ is a
multi-valued graph with gradient $\leq 1$.
\end{Cor}

Note that since $\Sigma$ is compact and embedded, the multi-valued
graphs given by Corollary \ref{c:conecor} spiral through the cone.
Namely, if a graph did close up, then the graph containing
$\Sigma_d$ would be forced to accumulate into it, contradicting
compactness.

 Another
result we need to apply Lemma \ref{l:regsing} is:

\begin{Lem}  \label{l:blowupnear}
See fig. 17.
There exists $c_0>0$ so:
Let $\Sigma_i\subset B_{R_i}$, $\partial \Sigma_i\subset \partial B_{R_i}$
be a sequence of embedded minimal disks with $R_i\to \infty$.
If $\Sigma_{d,i}\subset \Sigma_i$ is a sequence of $2$-valued
graphs over $D_{R_i/C}\setminus D_{\epsilon_i}$ with $\epsilon_i\to 0$
and $\Sigma_{d,i}\to \{x_3=0\}\setminus \{0\}$, then
\begin{equation}
\sup_{B_1 \cap \Sigma_i\cap \{x_3>c_0 \}}|A|^2\to \infty\, .
\end{equation}
\end{Lem}

\begin{figure}[htbp]
    \setlength{\captionindent}{20pt}
\begin{minipage}[t]{0.5\textwidth}
    \centering\input{uls17.pstex_t}
    \caption{Lemma \ref{l:blowupnear} - point with large curvature 
	in $\Sigma_i$ above the plane $x_3=c_0$ but near 
	the center of the $2$-valued graph $\Sigma_{d,i}$.}
    \end{minipage}\begin{minipage}[t]{0.5\textwidth}
    \centering\input{uls18.pstex_t}
    \caption{If Lemma \ref{l:blowupnear} failed, then by Corollary
\ref{c:conecor} the limit of the $\Sigma_i$'s
would contain a nonproper multi-valued graph contradicting
corollary 0.7 of \cite{CM7}.}
    \end{minipage}%
\end{figure}

\begin{proof}
Suppose not, see fig. 18; so assume that for each $c_0>0$, there is a sequence
of embedded minimal disks $\Sigma_i$ (and $C_1$ depending on both
$c_0$ and the sequence) with
\begin{equation}  \label{e:curvbd}
\sup_{B_1 \cap \Sigma_i\cap \{x_3>c_0 \}}|A|^2\leq C_1<\infty
\end{equation}
and $2$-valued graphs $\Sigma_{d,i} \subset \Sigma_i$ over
$D_{R_i/C}\setminus D_{\epsilon_i}$ with $\epsilon_i\to 0$,
$\Sigma_{d,i}\to \{x_3=0\}\setminus \{0\}$. Increasing
$\epsilon_i$ (yet still $\epsilon_i \to 0$) and replacing $R_i$ by
$S_i \to \infty$, we can assume $\Sigma_{d,i} \subset \{x_3^2 \leq
\epsilon_i^2\, (x_1^2+x_2^2)\}$ is  a $2$-valued graph over
$D_{4\e^{N_g} S_i}\setminus D_{\e^{-N_g} \epsilon_i/2}$ with
gradient $\leq \epsilon_i$ ($N_g$ given before \eqr{e:wantit}).

By Corollary \ref{c:conecor}, each component of $B_{2\e^{N_g} S_i}
\cap \Sigma_i\setminus ( \cone_{\delta_0}(0) \cup B_{\e^{-N_g}
\epsilon_i} )$ is a graph.  Hence, by the Harnack inequality,
 if $\alpha > 0$ is sufficiently small and
$q_i \in B_{S_i} \cap \Sigma_i \setminus ( \cone_{\alpha}(0) \cup
B_{2\epsilon_i} )$, then for $i$ large $\Sigma_i$ contains an
$(N_g+1)$-valued graph over $D_{\e^{N_g} |q_i|} \setminus
D_{\e^{-N_g} |q_i|/2}$ with gradient $\leq \epsilon < 1/(4\pi)$
and so $q_i$ is in the image of  $\{ |\theta| \leq \pi \}$ for
this graph. Consequently, each component of $B_{S_i} \cap \Sigma_i
\setminus ( \cone_{\alpha}(0) \cup B_{2\epsilon_i} )$ is a
multi-valued graph satisfying \eqr{e:wantit}.

Fix $h,\ell$ with $0< h < \alpha \, \ell$.  We get points $z_i \in
\{ x_3 = h, \, x_1^2 + y_1^2 = \ell^2 \} \cap \Sigma_i$ and
multi-valued graphs $z_i \in \Sigma_{1,i} \subset  \{ x_3
> 0 \} \cap \Sigma_i$ defined over
$S_{\ell/2,S_i/2}^{-3\pi, 3\pi + N_i}$, with $N_i \to \infty$, so
$z_i$ is in the image of $S_{\ell,\ell}^{-\pi , \pi}$, and so
$\Sigma_{1,i}$ spirals into $\{ x_3 = 0 \}$ (note that we have
assumed that it spirals down; we can argue similarly in the other
case).
 In particular, Theorem
\ref{t:main} applies, giving the other multi-valued graphs
$\Sigma_{2,i}$ so $\Sigma_{1,i}$ and $\Sigma_{2,i}$ spiral
together and so $\Sigma_{2,i}$ is the only part of $\Sigma_i$
between the sheets of $\Sigma_{1,i}$. Moreover, Theorem
\ref{t:main} also gives the short curves $\sigma_{\theta ,i}$
connecting these. It now follows from corollary $0.7$ of
\cite{CM7} that the separations of the graph $\Sigma_{1,i}$ at
$z_i$ go to $0$. Since this holds for all such $h$ and $\ell$, it
follows that $\Sigma_i\setminus \cone_{\alpha}(0)\to \cF$; where
$\cF$ is a foliation of $\RR^3\setminus \cone_{\alpha}(0)$ by
minimal annuli (all graphs over part of $\{x_3=0\}$).

 Theorem \ref{t:cm34}
gives $0<C_2<\infty$ so, given $r_0 > 0$, if $y_i\in
\Sigma_i\setminus B_{3r_0}$, $i$ is  large, and
\begin{equation}  \label{e:quadcm5}
|y_i|^2\,|A|^2(y_i)> C_2
\end{equation}
then there is a $2$-valued graph $\Sigma_{d,i}^{y_i}\subset
\Sigma_i\setminus B_{C_3 |y_i|}$   starting in $B_{C_4 |y_i|}(y_i)
\subset \{ x_3 > C_3 \, r_0 \}$ (by Theorem \ref{t:cm34},
$\Sigma_{d,i}^{y_i}$ starts in $B_{C_4 |y_i|}(y_i)$ where $C_4 =
C_4(C_2)$ and, by Corollary \ref{c:conecor}, $y_i \in
\cone_{\delta_0/2}(0)$). Let $C_2' = C_2'(C_2)>1$ be given by
Proposition \ref{c:0.4ofcm5} and set $r_0= 1/(4C_2')$.

Choose  $h_i , \ell_i \to 0$ with $\epsilon_i < \ell_i < r_0/4$,
$0< h_i < \alpha \, \ell_i$ and let $z_i , \Sigma_{1,i}
,\Sigma_{2,i}$ be as above. Since $\partial \Sigma_{i, z_i,2r_0}$
is a simple closed curve, it must pass between the sheets of
$\Sigma_{1,i}$. Since $\Sigma_{2,i}$ is the only part of
$\Sigma_i$ between the sheets of $\Sigma_{1,i}$, we can connect
$\Sigma_{1,i}$ and $\Sigma_{2,i}$ by  curves $\nu_i \subset
\partial \Sigma_{i, z_i,2 r_0}$ which are above $\Sigma_{1 , i}$.
We can now apply Proposition \ref{c:0.4ofcm5}
 to get the points
$y_i \in  B_{1/2}(z_i) \cap \Sigma_i \setminus B_{2r_0}(z_i)
\subset B_{1/2 + 4\ell_i} \setminus B_{3r_0}$  as in
\eqr{e:quadcm5}.

To get the desired contradiction, observe that if
$c_0<C_3 r_0$, then the $2$-valued graphs $\Sigma_{d,i}^{y_i}$
given by \eqr{e:quadcm5} and Theorem \ref{t:cm34}, have separation
$\geq C_5=C_5(C_1)>0$ (since $B_{C_4 |y_i|}(y_i) \subset \{ x_3 >
C_3 \, r_0 \}$). Namely, this separation is on a fixed scale
bounded away from zero even as $\Sigma_{d,i}^{y_i}$ extends out of
$\cone_{\alpha}(0)$,  contradicting $\Sigma_i\setminus
\cone_{\alpha}(0)\to \cF$ and the lemma follows.
\end{proof}

\section{Proof of Theorem \ref{t:t0.1}}

\begin{proof}
(of Theorem \ref{t:t0.1}). By Lemma
\ref{l:inftyornot}, after
passing to a subsequence (also denoted by $\Sigma_i$) we
can assume that for each $x\in \RR^3$ either
\begin{equation}  \label{e:curvbu}
\sup_{B_r(x)\cap \Sigma_i}|A|^2\to \infty\text{ for all $r>0$}\, ,
\end{equation}
or $\sup_i\sup_{B_r(x)\cap \Sigma_i}|A|^2<\infty$ for some $r>0$.
Let $\cS \subset \RR^3$ be the points where \eqr{e:curvbu} holds.
By assumption $B_1\cap\cS\ne \emptyset$.  So after a possible
translation we may assume that $0\in\cS$ and it follows easily
from the definition that $\cS$ is closed. By theorem $5.8$ of
\cite{CM4} (and Bernstein's theorem; see for instance theorem
$1.16$ of \cite{CM1}), there exists a subsequence $\Sigma_j$ and
$2$-valued graphs $\Sigma_{d,j}\subset \Sigma_j$ over
$D_{R_j/C}\setminus D_{\epsilon_j}$ with $\epsilon_j\to 0$ such
that $\Sigma_{d,j}\to \{x_3=0\}\setminus\{0\}$ (after possibly
rotating $\RR^3$). (This fixes the subsequence and the coordinate
system of $\RR^3$.)  Again by theorem $5.8$ of \cite{CM4} (and
Bernstein's theorem) for each $\cS_t\in \cS$ there are $2$-valued
graphs $\Sigma_{d,j}^t \subset \Sigma_j$ over
$D_{R_j/C}(\cS_t)\setminus D_{\epsilon_j}(\cS_t)$ with
$\epsilon_j\to 0$ such that $\Sigma_{d,j}^t \to
\{x_3=t\}\setminus\{\cS_t\}$.  Hence, by Corollary
\ref{c:conecor}, $\cS\subset \cone_{\delta}(\cS_t)$. By Lemma
\ref{l:blowupnear} (and scaling), for all $t\in x_3(\cS)$ and all
$\epsilon>0$, $\cS\cap \{t<x_3<t+\epsilon\}\ne \emptyset$,
$\cS\cap \{t-\epsilon<x_3<t\}\ne \emptyset$.  It follows from
Lemma \ref{l:regsing} that  $t\to \cS_t$ is a Lipschitz curve and
$\Sigma_j\setminus \cS\to \cF\setminus \cS$ in the
$C^{\alpha}$-topology for all $\alpha<1$ (and with uniformly
bounded curvatures on compact subsets of $\RR^3\setminus \cS$; see
also Appendix \ref{s:lamination}).
\end{proof}

\part{``The other half''}   \label{s:sotherhalf}

Theorem \ref{t:main}
will follow by first showing that if an embedded minimal disk contains a
multi-valued graph, then ``between the sheets'' of the graph the surface is
another
multi-valued graph - ``the other half''.   Second, we show an intrinsic
version of Theorem \ref{t:cm34} and, third, using this
intrinsic
version, we construct in Part \ref{p:shc} the short curves
connecting the two halves.

\section{``The other half'' of an embedded minimal disk}   \label{s:twohalves}

We show first that any point between
the sheets of a multi-valued graph must connect to it within a fixed
extrinsic ball:

\begin{Lem} \label{l:onecomp}
There exist $\epsilon_s > 0$, $C_s >2$ so: Let $0 \in \Sigma \subset
B_R$ be an embedded minimal disk with $\partial \Sigma \subset
\partial B_R$, $\Sigma_d \subset \{ x_3^2 \leq x_1^2 + x_2^2 \} \cap
\Sigma$ a $2$-valued graph over $D_{3r_0} \setminus D_{r_0}$ with
gradient $\leq \epsilon_s$. If $E_0$ is the region over $D_{2r_0}
\setminus D_{r_0}$ between the sheets  of $\Sigma_d$, then $E_0
\cap \Sigma \subset \Sigma_{0,C_s \, r_0}$.
\end{Lem}

\begin{proof}
Fix $\epsilon_s > 0$ small and $C_s$ large to be chosen.
 If the lemma fails,
then there are disjoint components $\Sigma_a , \Sigma_b$ of
$B_{C_s \, r_0} \cap \Sigma$ with $\Sigma_d \subset \Sigma_a$ and
$y \in E_0 \cap \Sigma_b$.
 By the maximum principle,
$\Sigma_a , \Sigma_b$ are disks. Let $\tilde \eta_y$ be the
vertical segment (i.e., parallel to the $x_3$-axis) through $y$
connecting the sheets of $\Sigma_d$. Fix a component $\eta_y$ of
$\tilde \eta_y \setminus \Sigma$ connecting $\Sigma_b$ to $\Sigma
\setminus \Sigma_b$. Let $\Omega$ be the component of $B_{C_s \,
r_0} \setminus \Sigma$ containing $\eta_y$ (so $\partial \Sigma_b$
and $\eta_y$ are linked in $\Omega$). \cite{MeYa} gives a stable
disk $\Gamma \subset \Omega$ with $\partial \Gamma = \partial
\Sigma_b$. Using the linking, $\Gamma$ intersects $\eta_y$.
\cite{Sc}, \cite{CM2} (cf. lemma I.0.9 of \cite{CM3}) give $C_s$
so any component $\Gamma_y$ of $B_{10r_0} \cap \Gamma$
intersecting $\eta_y$ is a graph with bounded gradient over some
plane; for $\epsilon_s$ small, this plane must be almost
horizontal.  Hence, $\Gamma_y$ is forced to ``cut the axis''
(i.e., intersect
 the curve in $\Sigma_d$ over $\partial D_{r_0}$ connecting
the top and bottom sheets), giving the desired contradiction.
\end{proof}

In the next proposition $\Sigma \subset B_{4R}$ with $\partial \Sigma \subset
\partial B_{4R}$
is an embedded minimal surface and $\Sigma_1 \subset \{ x_3^2 \leq
x_1^2 + x_2^2 \} \cap \Sigma$ an $(N+2)$-valued graph of $u_1$
over $D_{2R} \setminus D_{r_1}$ with  $|\nabla u_1| \leq
\epsilon$ and $N\geq 6$. Let $E_1$ be the region over $D_R \setminus D_{2 \,
r_1}$ between the top and bottom sheets of the concentric
$(N+1)$-valued subgraph in $\Sigma_1$.
 To be precise,
 $E_1$ is the set of all
$(r \cos \theta , r \sin \theta , t)$
with $2r_1 \leq r \leq R$, $(N-1) \pi \leq \theta < (N+1) \pi$, and
\begin{equation}    \label{e:defE1}
    u_1 (r , \theta ) < t < u_1 (r , \theta -2 N \pi)
     \, .
\end{equation}

\begin{Pro} \label{p:twosheets}
There exist $C_0 > C_s$ ,  $\epsilon_0 > 0$ so if $\Sigma$ is a
disk as above, $R \geq C_0 \, r_1$, and $\epsilon_0 \geq
\epsilon$, then $E_1 \cap \Sigma \setminus \Sigma_1$ is an
(oppositely oriented) $N$-valued graph $\Sigma_2$.
\end{Pro}

\begin{proof}
 Fix $z \in \Sigma_1$ over $\partial D_{r_1}$.
Since $\partial \Sigma_{z,2r_1}$ is a simple closed curve, it must
pass between the sheets of $\Sigma_{1}$ and hence through some
other component $\Sigma_2$ of $E_1 \cap \Sigma$.

The version of the ``estimate between the sheets'' given in
theorem III.2.4 of \cite{CM3} gives $\epsilon_0 > 0$ so that $E_1
\cap \Sigma$ is locally graphical (i.e., if $z \in E_1 \cap
\Sigma$, then $\langle \nn_{\Sigma}(z) , (0,0,1) \rangle \ne 0$).
It follows that each component of $E_1 \cap \Sigma$ is an
$N$-valued graph.

Fix a component $\Omega$ of $B_{4R} \setminus \Sigma$. We show
next that $\Sigma_2$ is the only other component of $E_1 \cap
\Sigma$ (i.e., $E_1 \cap \Sigma \subset \Sigma_1 \cup \Sigma_2$).
If not, then there is a third component $\Sigma_3$ which is also
an $N$-valued graph.  An easy argument (using orientations) shows
that there must then be a fourth component $\Sigma_4$ of $E_1 \cap
\Sigma$. Using that each $\Sigma_i$ is a multi-valued graph,
it follows easily that we can choose two of these four which cannot
be connected in $\Omega \cap E_1$; call these
 $\Sigma_{i_1} ,
\Sigma_{i_2}$.   The rest of this argument uses
 these components to find a stable
$\Gamma \subset \Omega$ which has points of large curvature
by Proposition \ref{c:0.4ofcm5}, contradicting the curvature
estimates from stability.   First, we
construct $\partial \Gamma$. Let $\sigma_j \subset
\Sigma_{i_j}$ be the images of $\{ \theta = 0 \}$ from $\{ x_1^2 +
x_2^2 = 4r_1^2 \}$ to $\partial B_R$ and set $y_j = \{ x_1^2 +
x_2^2 = 4r_1^2 \} \cap
\partial \sigma_j$. By Lemma \ref{l:onecomp}, $y_1
, y_2$ can be connected by a curve $\sigma_0 \subset B_{C_s \,
r_1} \cap \Sigma$. By the maximum principle, each component of
$B_R \cap \Sigma$ is a disk. Therefore, we can add a segment in
$\partial B_R \cap \Sigma$ to $\sigma_0 \cup \sigma_1 \cup
\sigma_2$ to get a closed curve $\sigma \subset \Sigma$.   A
result of \cite{MeYa} then gives a stable embedded minimal disk
$\Gamma \subset \Omega$ with $\partial \Gamma = \sigma$.

Now that we have $\Gamma$, we show that Proposition \ref{c:0.4ofcm5}
applies. Namely,
let (the disk) $\Gamma_{2C_s \,r_1}(\sigma_0)$ be the component of $B_{2C_s \,
r_1} \cap \Gamma$  containing $\sigma_0$, so that
$\partial \Gamma_{2C_s \,r_1}(\sigma_0)$ contains a curve
$\nu \subset \partial B_{2C_s \,r_1} $ connecting
$\sigma_1$ to $\sigma_2$.
Since $\sigma_1 , \sigma_2$ are in the middle sheets of
$\Sigma_{i_1} , \Sigma_{i_2}$ (and $\Gamma$ is stable),
$\Gamma$ contains two disjoint
$(N/2 - 1)$-valued graphs $\Gamma_1 , \Gamma_2$ in
$E_1$ which spiral together and $\nu$ connects these (note that
$E_1 \cap \Gamma$ may contain many components; at least
two of these, say $\Gamma_1 , \Gamma_2$, spiral together).
  Let
$\Gamma_0$ be the component of
$\Gamma_{R/2}(\sigma_0) \setminus (\nu \cup \Gamma_1 \cup \Gamma_2)$
which does not contain $\Gamma_{2C_s \,r_1}(\sigma_0)$. It is easy to
see that $\Gamma_0 \cap \partial \Gamma = \emptyset$; in fact,
if $x \in \Gamma_0$, then $\dist_{\Gamma} (x , \partial \Gamma) \geq |x|/2$.
Therefore,
  for $R/r_1$ sufficiently large,
Proposition \ref{c:0.4ofcm5} gives an interior point of large
curvature, contradicting the curvature estimate for stable
surfaces. We conclude that $E_1 \cap \Sigma \subset \Sigma_1 \cup
\Sigma_2$.  Finally, it follows easily  that $\Sigma_2$ is
oppositely oriented.
\end{proof}

The proof of Proposition
\ref{p:twosheets} simplifies when $\Sigma$ is in a slab.
In this case, \cite{Sc}, \cite{CM2} and the
gradient estimate (cf. lemma I.0.9 of \cite{CM3})
force $\Gamma$ to spiral indefinitely if it leaves $E_1$.

\section{An intrinsic version of Theorem \ref{t:cm34}}

We will first show a ``chord-arc'' type result
(relating extrinsic and intrinsic distances) assuming a curvature
bound on an intrinsic ball.

\begin{Lem}  \label{t:onesig}
(cf. lemma III.1.3 in \cite{CM5}). Given $R_0$, there exists $R_1$
so: If $0 \in \Sigma \subset B_{R_1}$ is an embedded minimal surface,
$\partial \Sigma \subset \partial B_{R_1}$, and $\sup_{ \cB_{R_1}
} |A|^2 \leq 4$, then $\Sigma_{0,R_0} \subset \cB_{R_1}$.
\end{Lem}

\begin{proof}
Let $\tilde{\Sigma}$ be the universal cover of $\Sigma$ and $\tpi
: \tilde{\Sigma} \to \Sigma$ the covering map.
With the definition of $\delta$-stable as in section $2$ of
\cite{CM4}, the argument of \cite{CM2}
(i.e., curvature estimates for $1/2$-stable surfaces) gives
 $C > 10$ so if
$\cB_{C R_0/2 }(\tilde{z}) \subset \tilde{\Sigma}$ is
$1/2$-stable
and  $\tpi(\tilde{z}) = z$, then $\tpi : \cB_{5 R_0}(\tilde{z}) \to
\cB_{5R_0}(z)$ is one-to-one and $\cB_{5 R_0}(z)$ is a graph with
$B_{4R_0}(z) \cap \partial \cB_{5R_0}(z) = \emptyset$.
Corollary 2.13 in \cite{CM4} gives $\epsilon = \epsilon (CR_0) > 0$
so if
$ |z_{1} - z_{2}| < \epsilon$ and $|A|^2 \leq 4$ on
(the disjoint balls)
$\cB_{ C R_0 }(z_i)$,
then each $\cB_{C R_0
/2}(\tilde{z}_{i}) \subset \tilde{\Sigma}$ is $1/2$-stable where
$\tpi (\tilde{z}_{i}) = z_i$.

We claim that there exists $n$ so $\Sigma_{0,R_0} \subset
\cB_{(2n+1) \, C R_0}$. Suppose not; we get a curve $\sigma
\subset \Sigma_{0,R_0} \subset B_{R_0}$ from $0$  to $\partial
\cB_{(2n+1)\,C R_0}$. For $i=1, \dots , n$, fix points $z_i \in
\partial \cB_{2i\,C R_0} \cap \sigma$.  It follows that the
intrinsic balls $\cB_{C R_0}(z_i)$ are disjoint, have centers in
$B_{R_0} \subset \RR^3$, and have $|A|^2 \leq 4$. In particular,
there exist $i_1 , i_2$ with $0< | z_{i_1} - z_{i_2}| < C' \,
R_0 \, n^{-1/3} < \epsilon$, and, by corollary 2.13 in
\cite{CM4},
each $\cB_{C R_0 /2}(\tilde{z}_{i_j}) \subset \tilde{\Sigma}$ is $1/2$-stable
where $\tpi ( \tilde{z}_{i_j}) = {z}_{i_j}$. By \cite{CM2}, each
$\cB_{5 R_0}(z_{i_j})$ is a graph with $B_{4R_0}(z_{i_j}) \cap
\partial \cB_{5R_0}({z}_{i_j}) = \emptyset$.  In particular,
$B_{R_0} \cap
\partial \cB_{5R_0}(z_{i_j}) = \emptyset$.
This contradicts that $\sigma \subset B_{R_0}$ connects
$z_{i_j}$ to $\partial \cB_{C R_0} (z_{i_j})$.
\end{proof}

An immediate consequence of
Lemma \ref{t:onesig}, is that we can improve Theorem
\ref{t:cm34} (and hence also, by an intrinsic blow-up
argument, Theorem \ref{t:blowupwinding0}) by observing that the
multi-valued graph can actually be chosen to be intrinsically
nearby where the curvature is large (as opposed to
extrinsically nearby):

\begin{Thm} \label{t:blowupwinding01}
Given $N\in \ZZ_+$, $\epsilon > 0$, there exist $C_1,\,C_2 , C_3>0$ so:
If $0\in \Sigma^2\subset B_{R}\subset \RR^3$ is an embedded
minimal disk, $\partial \Sigma\subset \partial B_{R}$, and
$\sup_{\cB_{r_0}}|A|^2\leq 4\,|A|^2(0)=4\,C_1^2\,r_0^{-2}$
for
some $0<r_0<R$, then there exists (after a rotation) an $N$-valued
graph $\Sigma_g \subset \Sigma\cap \{ x_3^2 \leq \epsilon^2 \,
(x_1^2 + x_2^2) \}$ over $D_{R/C_2} \setminus D_{r_0}$ with
gradient $\leq \epsilon$,
 separation $\geq C_3\,r_0$ over
$\partial D_{r_0}$,
and $\dist_{\Sigma}(0,\Sigma_g)\leq
2\,r_0$.
\end{Thm}

\begin{proof}
By combining theorems $0.4$ and $0.6$ of \cite{CM4}, we get
 $C_0 , C_2 , C_3$ so if
  $\sup_{\Sigma_{0,r_0}} |A|^2\leq 4\, |A|^2(0) = 4\,C_0^2 \, r_0^{-2}$,
then
we get
(after a rotation) an $N$-valued
graph $\Sigma_g \subset \Sigma\cap \{ x_3^2 \leq \epsilon^2 \,
(x_1^2 + x_2^2) \}$ over $D_{R/C_2} \setminus D_{r_0}$ with
gradient $\leq \epsilon$,
 separation $\geq C_3 \,r_0$ over
$\partial D_{r_0}$,
and which intersects
$\Sigma_{0,r_0}$.
Namely,
theorem $0.4$ of \cite{CM4} gives an initial $N$-valued graph
contained in $\Sigma_{0,r_0}$ and then theorem $0.6$ of \cite{CM4}
extends this out to $\partial D_{R/C_2}$.
Let $C_1$ be the
$R_1$
from Lemma \ref{t:onesig} with $R_0 = C_0$.
By rescaling, we can assume that $|A|^2(0) = 1$ and
$\sup_{\cB_{C_1}}|A|^2\leq 4$.  By Lemma \ref{t:onesig},
$\Sigma_{0,C_0} \subset \cB_{C_1}$, hence
$\sup_{\Sigma_{0,C_0}} |A|^2\leq 4$.
Theorems $0.4$ and $0.6$ of \cite{CM4} now
  give the desired $\Sigma_g$.
\end{proof}

A standard blowup argument gives points as in
Theorem \ref{t:blowupwinding01} (with $C_4 = 1$ and $s=r_0$):

\begin{Lem}     \label{l:l5.1}
(Lemma $5.1$ of \cite{CM4}).
Given $C_1 , C_4$, if $\cB_{C_1 C_4}(0) \subset \Sigma$
is an immersed surface and $|A|^2 (0) \geq 4$, then there exists
$\cB_{C_4 s}(z) \subset \cB_{C_1 C_4}(0)$ with
\begin{equation}    \label{e:defc1i}
    \sup_{\cB_{C_4 s}(z)}|A|^2\leq 4\,|A|^2(z)=4\,C_1^2\, s^{-2}
    \, .
\end{equation}
\end{Lem}

\begin{proof}
This follows as in Lemma $5.1$ of \cite{CM4}, except we define $F$
intrinsically on $\cB_{C_1 C_4}(0)$ by
$F = d^{2} \, |A|^2$
where
$d(x) = C_1 C_4 - \dist_{\Sigma} (x, 0)$
(so $F=0$ on $\partial \cB_{C_1 C_4}(0)$,
 $ F(0) \geq 4(C_1 C_4)^2$).  Let $F(z)$ be the maximum of $F$
and set $s = C_1 / |A|(z)$.  It follows that
$\sup_{\cB_{d(z)/2}(z)}|A|^2\leq 4\,|A|^2(z)$ and (using $F(z)
\geq (C_1 C_4)^2$) we get $2C_4 s \leq d(z)$, giving
\eqr{e:defc1i}.
\end{proof}

\part{The stacking and the proof of Theorem \ref{t:t2}} \label{p:p00}

This part deals  with how the
multi-valued graphs given by \cite{CM4} fit
together.    As mentioned
in the introduction, a general embedded minimal disk with large curvature at
some point should be
thought of as obtained by stacking such graphs on top of each other.

\section{Decomposing disks into multi-valued graphs}    \label{s:dd}

Fix $N > 6$ large, $1/10 > \epsilon > 0$ small. We will choose
$\epsilon_g > 0$  small depending on $\epsilon$
and then let $N_g = N_g (\epsilon_g)$  be
given by proposition II.2.12 of \cite{CM3}.
Below $\Sigma$ will be an embedded minimal disk.
Theorem
\ref{t:blowupwinding01} gives $C_1 , C_2 , C_3$
(depending on $\epsilon_g$, $N$, and $N_g$) so if
 $B_R(y)\cap \partial \Sigma = \emptyset$ and the pair
$(y,s)$ satisfies
\begin{equation}    \label{e:defc1ii}
    \sup_{\cB_{8 s}(y)}|A|^2\leq 4\,|A|^2(y)=4\,C_1^2\, s^{-2}
    \, ,
\end{equation}
then (after a rotation) we get
an $(N+ N_g +4)$-valued graph $\tilde{\Sigma}_1$ over
$D_{2\e^{N_g} R/C_2} (p)\setminus D_{\e^{-N_g}s/2}(p)$
with gradient $\leq \epsilon_g$,
separation $\geq C_3 \, s$ over $\partial D_{s}(p)$,
 and $\dist_{\Sigma} (y, \tilde{\Sigma}_1) \leq 2s$
(where $p = (y_1,y_2,0)$).
In particular, by proposition II.2.12 of \cite{CM3}
and  the version of the ``estimate between the sheets'' given in
theorem III.2.4 of \cite{CM3},
we can choose $\epsilon_g = \epsilon_g (\epsilon) > 0$ so that
(1) the concentric
$(N+3)$-valued subgraph $\hat{\Sigma}_1$ over
$D_{ R/C_2} (p)\setminus D_{s}(p)$ satisfies \eqr{e:wantit} and
(2) each component of $\Sigma$ between the sheets of
$\hat{\Sigma}_1$  (as in \eqr{e:defE1})
is an $(N+2)$-valued graph also satisfying
\eqr{e:wantit}.
In the remainder of this
section, $C_1 , C_2 , C_3$ will be fixed.

Let $\epsilon_0 , C_0$
  be from Proposition \ref{p:twosheets} and
suppose  $\epsilon < \epsilon_0$.
If $s < R/(8C_2 C_0)$ for
such a pair $(y,s)$, then
  Proposition \ref{p:twosheets} applies.
Let $\hat{E},E$ be the regions between the sheets of the concentric
$(N+2)$-valued and $(N+1)$-valued, respectively,
 subgraphs of $\hat{\Sigma}_1$ (over $D_{R/C_2}
(p)\setminus D_{s}(p)$).  By Proposition
\ref{p:twosheets} (and (2) above),
$\hat{E} \cap \Sigma \setminus \hat{\Sigma}_1$
is an $(N+1)$-valued graph $\hat{\Sigma}_2$;
similarly,
$E \cap \Sigma \setminus \hat{\Sigma}_1$ is an $N$-valued
graph $\Sigma_2 \subset \hat{\Sigma}_2$.
Let $\Sigma_1 \subset \hat{\Sigma}_1$ be
the concentric $N$-valued subgraph.
Since $\partial
\Sigma_{y,4s}$ is a simple closed curve, it must pass through $E
\setminus \Sigma_1$. Therefore, since   $\Sigma_{2}$ is the only
other part of $\Sigma$ in $E$, we can connect $\Sigma_{1}$ and
$\Sigma_{2}$ by  curves $\nu_{\pm} \subset \partial B_{4s}(y) \cap
\Sigma$ which are above and below $E$, respectively. This gives
components $\Sigma_{\pm}$ of $\Sigma_{y,R/(2C_2)} \setminus
(\Sigma_1 \cup \Sigma_2 \cup \nu_{\pm})$ which do not contain
$\Sigma_{y,s}$ and which are above  and below $E$, respectively
(these will be the $\Sigma_0$'s for Proposition \ref{c:0.4ofcm5}).

Given a pair satisfying \eqr{e:defc1ii}, Proposition
\ref{c:0.4ofcm5} and
Lemma \ref{l:l5.1} easily
give two nearby pairs
(one above and one below):

\begin{Lem} \label{l:dbl}
There exists $\cdbl > 1$ so: If
 $0 \in \Sigma \subset B_{3R}$ is an embedded minimal disk with
$\partial \Sigma \subset \partial B_{3R}$, $(0,s)$ satisfies
\eqr{e:defc1ii}, and $s < \min \{ R /(2 \cdbl) , R/(8 C_2 C_0) \}$,
 then we get $(y_{-},s_{-})$ also satisfying \eqr{e:defc1ii}
with  $y_{-} \in \Sigma_{-}$ and
 $\Sigma_{y_- , 4s_-} \subset \Sigma_{0,\cdbl s} \setminus B_{4 s}$.
Moreover, the $N$-valued graphs
 corresponding to $(0,s)$, $(y_{-},s_{-})$ are  disjoint.
\end{Lem}

\begin{proof}
Proposition
\ref{c:0.4ofcm5} gives $\cdbl = \cdbl (C_1)$ and
$z \in \Sigma_{0,\cdbl s/2} \cap \Sigma_{-}
\setminus B_{8s}$ with $|z|^2 \, |A|^2(z) \geq 4 (8C_1)^2$.
Since $\hat{E} \cap \Sigma$ consists of the multi-valued graphs
$\hat{\Sigma}_1 , \hat{\Sigma}_2$, we have $|x|^2 \, |A|^2(x) \leq C$
 on $\hat{E} \cap B_{\cdbl s} \cap \Sigma_- \setminus B_{2s}$ for
$C$ small ($C$ can be made
arbitrarily small by choosing $\epsilon$ even smaller).
Hence, $z \notin \hat{E}$ and so
$\cB_{|z|/2}(z)  \cap E = \emptyset$.
Applying Lemma \ref{l:l5.1} on $\cB_{|z|/2}(z)$, we get
 $(y_- , s_-)$ satisfying \eqr{e:defc1ii} with
$\cB_{8s_-}(y_-) \subset \cB_{|z|/2}(z)$ ($\subset \Sigma_- \setminus E$).
It follows that $\Sigma_{y_- , 4s_-} \subset \Sigma_{0,\cdbl s}\setminus B_{4 s}$ and
the corresponding $N$-valued graphs are disjoint.
\end{proof}

Let $\cdbl $ be given by Lemma \ref{l:dbl}. Iterating the
construction of Lemma \ref{l:dbl}, we can decompose an embedded
minimal disk into basic building blocks ordered by heights (the
points $p_i$ in Corollary \ref{c:decompo}  are the projections
 to $\{ x_3 = 0\}$  of the blowup points $y_i$):

\begin{Cor} \label{c:decompo}
There exist $C_5 > 1$, $\tilde{C}_3 > 0$ so: Let $\Sigma \subset
B_{C_5 R}$ be an embedded minimal disk, $\partial \Sigma \subset
\partial B_{C_5 R}$. If $(y_0,s_0)$ satisfies
\eqr{e:defc1ii} with $B_{\cdbl}(y_0) \subset B_{R}$,
then there exist $\{ (y_i, s_i) \}$
(for $i>0$) satisfying
\eqr{e:defc1ii} with $y_i \in \Sigma_-$
and corresponding (disjoint)
 $N$-valued graphs $\Sigma_i \subset \Sigma$
of $u_i$ over $D_{2R}(0) \setminus D_{2s_i}(p_i)$
with
gradient $\leq 2 \epsilon$, separation $\geq \tilde{C}_3  \, s_i$ over
$\partial D_{2s_i}(p_i)$,
\begin{alignat}{5}
\text{if $i< j$ and both $u_i, u_j$ are defined }
\text{at $p$, then } u_j (p) < u_i (p)
\, , \label{e:obh} \\
\Sigma_{y_{i+1} , 4 s_{i+1}} \subset \Sigma_{y_i , \cdbl s_i}
\setminus B_{4 s_i}(y_i) \text{ and }
  \cup_i B_{\cdbl s_i}(y_i) \setminus B_{R} \ne \emptyset\, .
\label{e:leaves}
\end{alignat}
\end{Cor}

\begin{proof}
Starting with $(y_0, s_0)$, we can apply Lemma
\ref{l:dbl}  repeatedly, until the second part of
\eqr{e:leaves} holds,
to find bottom $N$-valued graphs
giving \eqr{e:obh} and the first part of \eqr{e:leaves}.
 Each $N$-valued
graph is a graph over some plane with gradient $\leq \epsilon$.
Since $\Sigma$ is embedded, we can take these to be graphs over
a fixed plane with gradient $\leq 2 \epsilon$ (after possibly taking
$C_5 > 3C_2+1$ larger).  $\tilde{C}_3  > 0$ is now just a
fixed fraction of $C_3$.
\end{proof}

In the next lemma and  corollary,  $\Sigma \subset
B_{C_5 R}$ is an embedded minimal disk, $\partial \Sigma \subset
\partial B_{C_5 R}$.

\begin{Lem}     \label{l:mis}
If $(y,s)$ satisfies \eqr{e:defc1ii}, $B_s (y) \subset
B_{R/2}$, then the corresponding
$2$-valued graph over
$D_{R}(0)\setminus D_s(p)$ (after a rotation)
 has separation $\geq C_3  \, (s/R)^{\epsilon}  \, s/2$ at
$\partial D_{R}(0)$.
\end{Lem}

\begin{proof}
By the discussion around \eqr{e:defc1ii}, the
separation $|w|$ is  $\geq {C}_3   \, s$ at $\partial D_{s}(p)$ and
$|\nabla \log |w|| \leq \epsilon / \rho_p$ on $D_{2R}(p) \setminus
D_s (p)$.  Since $D_s(p) \subset D_{R/2} (0)$, integrating gives
\begin{equation}        \label{e:intgdt}
    \min_{\partial D_{R}(0)}  |w| \geq \min_{D_{2R}(p)\setminus
    D_{R/2}(p)} |w| \geq {C}_3   \, (s/(2R))^{\epsilon} \, s
        \, .
\end{equation}
\end{proof}

\begin{Cor}     \label{c:mis}
There exists $\cmis>0$ so if $(0,s)$ satisfies \eqr{e:defc1ii} and
$\sup_{B_{\ell}(0) \cap \Sigma_-}|A|^2\leq 5C_1^2 s^{-2}$
for some $ 4 \cdbl^2 \, s<{\ell}<R$, then there exists
$(z,r)$ satisfying \eqr{e:defc1ii} with $\Sigma_{z,r} \subset \Sigma_{0,\ell/2}$,
so the separation at $\partial D_{\ell}(0)$ between the $2$-valued graphs
$\Sigma_0$, $\Sigma_z$, corresponding to $(0,s)$, $(z,r)$, is $\geq \cmis
\,    (s/{\ell})^{\epsilon}  \, {\ell}$, and $\Sigma_z \subset \Sigma_-$.
\end{Cor}

\begin{proof}
Set $(y_0,s_0)=(0,s)$ and let $(y_i,s_i)$, $\Sigma_i$, $u_i$, $p_i$
be given by Corollary \ref{c:decompo}.  Let $i_0$ be the first $i$ with
$B_{\cdbl s_{i_0}} (y_{i_0}) \setminus B_{{\ell}/2}(0) \ne \emptyset$.
Set $(z,r) = (y_{i_0-1} , s_{i_0-1})$.
It follows  for $i < i_0$ that $B_{s_i}(y_i) \subset B_{{\ell}/2}(0)$
and
$s_i \geq s/2$ since $\sup_{B_{\ell}(0)\cap \Sigma_-}|A|^2\leq 5C_1^2 s^{-2}$.
Hence, by  Lemma \ref{l:mis} (as in Corollary \ref{c:decompo}),
$\Sigma_i$ has separation
$\geq \tilde{C}_3   \, (s/\ell)^{\epsilon} \, s_i / 4$ at
$\partial D_{{\ell}}(0)$
for $i < i_0$.   By \eqr{e:leaves},
${\ell}/4 \leq  \sum_{i \leq i_0} C_4 s_i
\leq (1+C_4) \sum_{i < i_0} C_4 s_i $.
Since the $\Sigma_i$'s are ordered by height,
the separation at $\partial D_{\ell}(0)$ between
$\Sigma_0$ and $\Sigma_z=\Sigma_{i_0-1}$ is
 $\geq \sum_{i < i_0} \tilde{C}_3   \, (s/\ell)^{\epsilon} \, s_i  /4
\geq \cmis \,
    (s/{\ell})^{\epsilon}  \, {\ell}$.
\end{proof}

\section{Stacking  multi-valued graphs and Theorem  \ref{t:t2}}
\label{s:somvg}

 If $(y,s)$ satisfies \eqr{e:defc1ii}, then $\Sigma_y$ is
the corresponding $2$-valued graph and $\Sigma_{y,-}$ the portion of
$\Sigma$ below $\Sigma_y$.  Given $C>8$, we will consider such pairs which
in addition satisfy
\begin{equation}   \label{e:pairs}
\sup_{B_{Cs}(y) \cap \Sigma_{y,-}}|A|^2 \leq 4\,|A|^2(y)=4C_1^2
s^{-2}\, .
\end{equation}

Using Section \ref{s:dd},
we show next that a pair $(0,s)$ satisfying \eqr{e:pairs}
has a nearby pair with a definite height
 below $\Sigma_0$.
$\Sigma\subset B_{C_5 R}$, $\partial \Sigma \subset
\partial B_{C_5 R}$ is an embedded minimal disk.

\begin{Pro}   \label{p:lift}
See fig. 9.
There exist $C , \bar{C} > 10 \, \cdbl^2$ and $\delta > 0$ so if $(0,s)$
satisfies \eqr{e:pairs} with $s< R/\bar{C}$,
 $\Sigma_0\subset \Sigma$ is  over $D_R\setminus D_s$ (without a
rotation), and $\nabla u( (Rs)^{1/2} , 0) = 0$, then
we get $(y,t)$ satisfying \eqr{e:pairs} with $y\in \cone_{\delta}(0)
\cap \Sigma_- \setminus B_{Cs/2}$.
\end{Pro}

\begin{proof}
We will choose $C$
large below and then set $\delta = \delta (C) > 0$, $\bar{C} = \bar{C}(C)$.
Note first  that (since $\nabla u( (Rs)^{1/2} , 0) = 0$),
 corollary $1.14$ of
\cite{CM7} gives that  $|\nabla u (\rho , \theta)| \leq  C_a \,
  (\rho / s)^{-5/12}$ for
$s
\leq \rho \leq (Rs)^{1/2}$.  Integrating this, we get for $s
\leq \rho \leq (Rs)^{1/2}$
\begin{equation}    \label{e:sigmadincone0}
      |u(\rho , \theta)| \leq s +
        C_a \,
 \int_s^{\rho}
     (\tau / s)^{-5/12}  \, d \tau
\leq (1+2C_a )\,(s/\rho)^{5/12}\,\rho\, .
\end{equation}

Proposition
\ref{c:0.4ofcm5} gives $\ccmf  (C_1 , C)$,
$z_0 \in B_{\ccmf  s} \cap \Sigma_{-}
\setminus B_{4s}$ with $ |A|^2(z_0) \geq 5 \, C^2 \,  C_1^2 \, |z_0|^{-2}$.
 Set
 \begin{equation}   \label{e:defx}
    \cA=\{x\in B_{\ccmf  s} \cap \Sigma_- \,|\,|A|^2(x)
\geq 5 \, C^2 \, C_1^2 |x|^{-2}   \}\, ,
\end{equation}
(so $z_0 \in \cA$)
and let $x_0\in \cA$ satisfy $|x_0|=\inf_{x\in\cA}|x|$. So
$|A|^2 \leq 5
\,  C_1^2 \, s^{-2}$
    on $B_{|x_0|} \cap \Sigma_-$ by \eqr{e:pairs} and
$Cs \leq |x_0|\leq \ccmf  \,s $.  An obvious extrinsic version of
  Lemma \ref{l:l5.1} (cf. Theorem \ref{t:t2.1}) gives
      $(y,t)$ satisfying \eqr{e:pairs} with
$B_{Ct} (y) \subset B_{|x_0|/2} (x_0)$.
We can assume $|p|\geq 4|y|/5$.

Since $|A|^2 \leq 5 \, C_1^2 \, s^{-2}$ on
$B_{|x_0|/2} \cap \Sigma_-$ and $(0,s)$ satisfies \eqr{e:pairs}
hence \eqr{e:defc1ii}, Corollary \ref{c:mis} (with $\ell = |p|$)
 gives $(z,r)$ also satisfying
\eqr{e:defc1ii} with $\Sigma_z \subset \Sigma_{-}$,
$\Sigma_{z,r} \subset \Sigma_{0,|p|/2}$,
and so the separation between  $\Sigma_{0}$
and $\Sigma_z$ is at least $C_c
\,    (s/|y|)^{\epsilon} \, |y|$ at $p$.
However, since $\Sigma$ is embedded, then $\Sigma_y$ must be
below both $\Sigma_0$ and $\Sigma_z$.  Combining this with
\eqr{e:sigmadincone0} gives
\begin{equation}    \label{e:gainb}
   |x_3| (y)/|y| \geq  C_c\,(s/|y|)^{\epsilon}
-(1+2\,C_a )(s/|y|)^{5/12}\,\, .
\end{equation}
 Since $C\leq 2|y|/s\leq 3\ccmf $, by choosing $C$
sufficiently large and then setting
$\bar{C} = \bar{C} (C,C_5)$, $\ccmf  = \ccmf  (C)$ (where
$\bar{C}$ is chosen so $\ccmf  \,s \leq (Rs)^{1/2}$) the proposition follows from
\eqr{e:gainb}.
\end{proof}

We show next Theorem \ref{t:t2}. Namely, iterating
Proposition \ref{p:lift}, we show that if the curvature
of an embedded minimal disk were large at a point, then it
would be forced to
grow out of the half-space it was assumed to lie in.
First we need:

\begin{Lem}   \label{l:apriori}
Given $C$, $\delta>0$, there exists $\epsilon_1>0$ so: Let
$\Sigma\subset B_{2r_0}\cap \{x_3>0\}$ be an embedded minimal
disk, $\partial \Sigma\subset \partial B_{2r_0}$, and
$\sup_{\Sigma \cap \{x_3\leq \delta\,r_0\}}|A|^2\leq C\,r_0^{-2}$,
then  $\sup_{\Sigma'}|A|^2\leq r_0^{-2}$ for all components
$\Sigma'$ of $B_{r_0}\cap \Sigma$ which intersect $B_{\epsilon_1
r_0}$.
\end{Lem}

\begin{proof}
If $y \in B_{r_0}\cap \Sigma \cap \{x_3\leq \delta\,r_0 /
4\}$, then
$\sup_{\Sigma_{y,\delta r_0/2}} |\nabla x_3|^2 \leq C
\, x_3^2 (y) \, \delta^{-2} \, r_0^{-2}$
(by the gradient estimate) and hence
$\Sigma_{y,\delta r_0/2}$ is a graph for  $x_3 (y) /( \delta \,
r_0)$ small; cf. Lemma \ref{l:l2.1}.  The lemma follows by
applying this to a chain of balls as in  the proof of lemma $2.10$
in \cite{CM8} or the theorem in \cite{CM10}.
\end{proof}

Let $C_1, \dots , C_6$ be as above and $\delta , C , \bar{C}$ be from
Proposition \ref{p:lift}.

\begin{proof}
(of Theorem \ref{t:t2})
By Lemma \ref{l:apriori} (and scaling),
it suffices to find $d>0$
and $\hat{C} > 1$
so  if $\Sigma \subset B_{4 C_5 \hat{C} R} \cap \{ x_3
> 0 \}$ and $\partial \Sigma \subset
\partial B_{4 C_5 \hat{C} R}$, then
\begin{equation}    \label{e:stef}
\sup_{B_{dR}\cap \Sigma}|A|^2 \leq 4 \, C^2 \, C_1^2 \, (dR)^{-2}\, .
\end{equation}

Suppose \eqr{e:stef} fails; we will get a contradiction.
An obvious extrinsic version of
  Lemma \ref{l:l5.1} gives
  $(y_0,s_0)$ satisfying \eqr{e:pairs}
with $B_{C s_0} (y_0) \subset B_{2dR}$.  Let $\Sigma_0$ be the corresponding
$N$-valued graph of $u_0$ over $D_{\hat{C} R} \setminus D_{s_0}(p_0)$
and $\Sigma_-$ the portion of $\Sigma$ below $\Sigma_0$.

To apply Proposition \ref{p:lift} we will need that if $(y,s)$
satisfies \eqr{e:pairs} with $y \in B_{2R} \cap \Sigma_-$ (where
$\Sigma_y$ is a graph of $u$ over $D_{\hat{C} R} \setminus
D_{s}(p)$), then
\begin{equation}  \label{e:edge}
s\leq R/ \bar{C}  {\text{ and }} |\nabla u ( (\hat{C}Rs)^{1/2} ,
0)| < \delta / 4\, .
\end{equation}
To see \eqr{e:edge}, observe first that
by proposition
II.2.12 of \cite{CM3} (since $u_0 > 0$),
$\sup_{D_{6R}} u_0 \leq  2dR (6/d)^{\epsilon} \leq 12 \, R \, d^{1-\epsilon}$.
It follows that
$B_{6R} \cap \Sigma_- \subset \{ 0 < x_3 <
12 \, R \, d^{1-\epsilon} \}$; hence, $s\leq C_a \, R \, d^{1-\epsilon}$ and
by the gradient estimate
$\sup_{\partial D_{4R}} |\nabla u| \leq C_b \, d^{1-\epsilon}$.
Lemma $1.8$ of \cite{CM7} and the mean value inequality (as in
corollary $1.14$ of
\cite{CM7}) gives  $|\Hess_u| \leq C_c (\hat{C} R)^{-5/12} \, \rho^{-7/12}$ for
$(\hat{C} R s)^{1/2} \leq \rho \leq \hat{C} R$.  Combining these gives at
$\rho = (\hat{C} Rs )^{1/2}$, $\theta = 0$
\begin{equation}    \label{e:graduibd}
    |\nabla u | \leq C_b \, d^{1-\epsilon} +
    C_c  (\hat{C} R)^{-5/12} \, \int_{0}^{8R}
     t^{-7/12} \, dt =  C_b \, d^{1-\epsilon} + C_c' \hat{C}^{-5/12} \, .
\end{equation}
In particular, for $d>0$ small and $\hat{C}$ large, \eqr{e:graduibd} gives
\eqr{e:edge}.

Repeatedly applying Proposition \ref{p:lift} (using \eqr{e:edge}) gives
$(y_{i+1}, s_{i+1})$ satisfying  \eqr{e:pairs} with
$y_{i+1} \in \cone_{\delta/2}(y_i)
\cap \Sigma_- \setminus B_{Cs_i /2}(y_i)$.
After choosing $d>0$ even smaller,
it follows that
the $y_i$'s must leave the half-space before they leave $B_R$.
\end{proof}

\begin{proof}
(of Corollary \ref{c:barrier}).  Using $\Sigma_1 \cup \Sigma_2$ as
a barrier, \cite{MeYa}, and a linking argument (cf. Lemma
\ref{l:onecomp}) give
 a stable surface $\Gamma \subset
B_{cr_0} \setminus (\Sigma_1 \cup\Sigma_2)$ with
$\partial \Gamma \subset \partial
B_{cr_0}$ and $B_{\epsilon r_0} \cap
\Gamma \ne \emptyset$.  Estimates for stable surfaces give a graphical
component of $B_{2r_0} \cap \Gamma$ which intersects $B_{\epsilon r_0}$.
The corollary now follows from Theorem \ref{t:t2}.
\end{proof}

\part{The short connecting curves and Theorem \ref{t:main}} \label{p:shc}

We first combine Lemmas \ref{l:onecomp} and \ref{t:onesig} to
see that any curve in $\Sigma$ which intersects both above and
below a multi-valued graph (with a curvature bound on an intrinsic
ball) connects to it in a fixed intrinsic ball:

\begin{Cor}   \label{c:onesig}
Given $C_1$, there exists $C_4 > 1$ so: Let $\Sigma , \Sigma_d,
E_0 , r_0$ be as in Lemma \ref{l:onecomp}.  If $\sup_{\cB_{C_4
r_0}}|A|^2 \leq 4 \, C_1^2 \, r_0^{-2}$ and $\eta \subset B_{2r_0}
\cap \Sigma$ connects points in $\partial B_{2r_0}$ above and
below $E_0$, then $\eta \subset \cB_{C_4 r_0}$.
\end{Cor}

\begin{proof}
Let $\Sigma_{2r_0}(\eta)$ be the component of $B_{2r_0} \cap
\Sigma$ containing $\eta$.  By the maximum principle,
$\Sigma_{2r_0}(\eta)$ is a disk and so $\partial
\Sigma_{2r_0}(\eta)$ must pass through $E_0$ (to connect the
points on opposite sides of $E_0$). Hence, by Lemma
\ref{l:onecomp}, $\Sigma_{2r_0}(\eta) \subset \Sigma_{0,C_s
 r_0}$.   Finally, by Lemma \ref{t:onesig},
$\Sigma_{0,C_s
 r_0} \subset  \cB_{C_4 r_0}$, giving the corollary.
\end{proof}

\begin{figure}[htbp]
    \setlength{\captionindent}{20pt}
    \begin{minipage}[t]{0.5\textwidth}
    \centering\input{uls19.pstex_t}
    \caption{If Theorem \ref{t:main} fails, then there are points
    $y_1\in \Sigma_1$ and $y_2\in \Sigma_2$ in consecutive sheets
    which are intrinsically far apart.}
    \end{minipage}\begin{minipage}[t]{0.5\textwidth}
    \centering\input{uls20.pstex_t}
    \caption{Proof of Theorem \ref{t:main}.  The blowup points $z_1$,
    $z_2$ and the corresponding multi-valued graphs $\hat \Sigma_1$,
    $\hat \Sigma_2$ and the curves $\eta_i$ connecting
    $y_i$ with $\hat \Sigma_i$.}
    \end{minipage}

\end{figure}

\begin{proof}
(of Theorem \ref{t:main}). Fix $\epsilon >0$ with
$\epsilon < \min \{ \epsilon_0 ,
\epsilon_s \}$ ($\epsilon_0$ given by Proposition \ref{p:twosheets}
and $\epsilon_s$ from
Lemma \ref{l:onecomp}).  Choose $N_0 , R_0$ large so that
Proposition \ref{p:twosheets} gives ``the other half'' $\Sigma_2$.
If $\Sigma_1$ comes from an intrinsic blow up point, then it
follows from Lemma \ref{l:onecomp},  that there are short curves
connecting $\Sigma_1$ and $\Sigma_2$. While it is a priori not clear that
every multi-valued graph arises this way, Theorem
\ref{t:blowupwinding01} implies that every multi-valued graph is
intrinsically near one of these. We use this below to produce the
short curves $\sigma_{\theta}$ in general.

Suppose that no $\sigma_{\theta}$ with length $\leq C$ exists for
some $\theta$; we get $y_i\in \{ x_1^2 + x_2^2 = 1 \}\cap
\Sigma_i$ for $i=1,2$ with $\dist_{\Sigma} (y_1,y_2) > C$ and so
$y_1 , y_2$
 are in consecutive sheets of $\Sigma$ (i.e., $y_1$
and $y_2$ can be connected by a segment parallel to the $x_3$-axis
which does not otherwise intersect $\Sigma$). See fig. 19.  We will get a
contradiction from this  for $C$ large.

Since $\partial \Sigma_{y_1,4}$ is a simple closed curve, it must
pass through $E \setminus \Sigma_1$.  See fig. 20.
Therefore, since $\Sigma_{2}$
is the only other part of $\Sigma$ in $E$, we can connect
$\Sigma_{1}$ and $\Sigma_{2}$ by a curve in $\partial
\Sigma_{y_1,4}$.  Connecting the endpoints of this curve to $y_1$
and $y_2$ gives  a curve $\eta \subset \Sigma_{y_1,4}$ from $y_1$
to $y_2$. Since  $\cB_{4}(y_i)$ is not a graph, $\sup_{\cB_4
(y_i)}|A|^2 \geq C_0 > 0$. Let $C_1 , C_2$
(depending on $\epsilon$ and some fixed $N > 6$)
be given by  Theorem
\ref{t:blowupwinding01} and
 $C_4 = C_4 (C_1)$ given by Corollary  \ref{c:onesig}.
Lemma \ref{l:l5.1} gives pairs $(z_i,s_i)$ satisfying
\eqr{e:defc1i}
 with $\cB_{C_4 s_i}(z_i) \subset
\cB_{C'}(y_i)$ where $C'$ does not depend on $C$. Let
$\hat{\Sigma}_{1} , \hat{\Sigma}_{2}$  be the  multi-valued graphs
given by Theorem \ref{t:blowupwinding01} and
 $\hat{E}_i$ the regions between
the sheets.  Since
    $\dist_{\Sigma}(z_i, \hat{\Sigma}_i )\leq 2\,s_i$, we can
choose  curves $\eta_i$ from $y_i$ to $\cB_{2s_i}(z_i) \cap
\hat{\Sigma}_i$ with length $\leq C'$. Combining Corollary
\ref{c:onesig}, ${\text{Length}}(\eta_i)\leq C'$, and
$\dist_{\Sigma} (y_1,y_2) > C$, it is easy to see
 that, for $C$ large, $\eta_1$ intersects only one side
of $\hat{E}_2 \cup B_{2s_2}(z_2)$;
similarly, $\eta_2$ intersects only one side of $\hat{E}_1 \cup
B_{2s_1}(z_1)$.

We will next find a third pair $(z_3 , s_3)$ satisfying
\eqr{e:defc1i} which is between $\hat{E}_1 \cup B_{2s_1}(z_1)$ and
$\hat{E}_2 \cup B_{2s_2}(z_2)$ but which is intrinsically far from
$\eta_1 , \eta_2$; Corollary \ref{c:onesig} will then give a
contradiction.  Since $\dist_{\Sigma} ( \eta_1 ,
\eta_2 ) > C-2C'$, $\eta_1$ intersects only one side of $\hat{E}_2
\cup B_{2s_2}(z_2)$,  $\eta_2$ intersects only one side of
$\hat{E}_1 \cup B_{2s_1}(z_1)$, and $\eta \subset \Sigma_{y_1,4}$
connects $y_1 , y_2$,  it is easy to see that
there is a point $y_3 \in
\Sigma_{y_1,4}$ with $\dist_{\Sigma} (y_3 , \, \{\eta_1 , \eta_2
\} ) > (C-2C') / 2$ and so $y_3$ is between $\hat{E}_1 \cup
B_{2s_1}(z_1)$ and $\hat{E}_2 \cup B_{2s_2}(z_2)$.  (This last
condition means that there is a curve $\eta_{y_3}$ from
$B_{2s_1}(z_1)$ to $B_{2s_2}(z_2)$ so $y_3 \in \eta_{y_3}$ and
$\eta_{y_3}$ intersects only one side of each of $\hat{E}_1 \cup
B_{2s_1}(z_1)$, $\hat{E}_2 \cup B_{2s_2}(z_2)$.)  As above, Lemma
\ref{l:l5.1} gives a pair $(z_3 , s_3)$ satisfying \eqr{e:defc1i}
with   $\cB_{C_4 s_3} (z_3) \subset \cB_{C'}(y_3)$ and then
Theorem \ref{t:blowupwinding01} gives corresponding
$\hat{\Sigma}_{3}$, $\hat{E}_3$. Since $C'$ does not depend on
$C$, we can assume that
\begin{equation}    \label{e:3far}
    \dist_{\Sigma} (z_3 , \, \{\eta_1 , \eta_2 \} ) > C / 4 \, .
\end{equation}
It follows easily from Corollary \ref{c:onesig} that $\hat{E}_3$
is between $\hat{E}_1$ and $\hat{E}_2$ (since $\hat{\Sigma}_3$ is
close to $y_3$ and $y_3$ is far from $\hat{\Sigma}_1 ,
\hat{\Sigma}_2$). Moreover, it is easy to see that at least one of
$\eta_1 , \eta_2$ must intersect both sides of $\hat{E}_3 \cup
B_{2s_3}(z_3)$ and, therefore, Corollary \ref{c:onesig} gives
\begin{equation}    \label{e:3notfar}
    \dist_{\Sigma} (\hat{\Sigma}_3 , \, \{\eta_1 , \eta_2 \} ) \leq C''
\end{equation}
($C''$ independent of $C$).  For $C$  large, \eqr{e:3far}
contradicts \eqr{e:3notfar}, giving the theorem.
\end{proof}

\appendix

\section{One-sided Reifenberg condition and curvature estimates}
\label{s:ref}

We will show here curvature estimates for minimal
hyper-surfaces, $\Sigma^{n-1}\subset M^n$, which on all
sufficiently small scales lie on one side of, but come close to,
a hyper-surface with small curvature. Such a
minimal hyper-surface is said to satisfy the {\it
one-sided Reifenberg condition}. Note that no assumption on the
topology is made.
Inspired by the classical Reifenberg condition
(cf. \cite{ChC} and references therein)
we make the definition:

\begin{Def} \label{d:dosr}
A subset, $\Gamma$, of $M^n$  satisfies the
$(\delta,r_0)$-{\it one-sided Reifenberg condition} at $x\in
\Gamma$ if for every $0<\sigma\leq r_0$ and every $y\in
B_{r_0-\sigma}(x)\cap \Gamma$, there is a connected
hyper-surface, $L_{y,\sigma}^{n-1}$, with $\partial
L_{y,\sigma}\subset \partial B_{\sigma}(y)$,
\begin{equation}
B_{\delta \sigma}(y)\cap L_{y,\sigma}\ne \emptyset , \,
\sup_{B_{\sigma}(y)\cap L} |A_L|^2 \leq \delta^2 \, \sigma^{-2}\, ,
\end{equation}
and the component of $B_{\sigma}(y)\cap \overline{\Gamma}$ through
$y$ lies on one side of $L_{y,\sigma}$.
\end{Def}

\begin{Lem} \label{l:l2.1}
There exist $\rc (i_0,k,n) >0$, $0 < \epsilon_0 < 1$, $C=C(n)$ so for
$\epsilon \leq  \epsilon_0$, $r_0 \leq \rc$: If $z\in
\Sigma^{n-1}\subset B_{r_0}=B_{r_0}(z)\subset M^n$ is an embedded
minimal hyper-surface, $\partial \Sigma \subset \partial
B_{r_0}$, and there is a connected hyper-surface, $L^{n-1}$,
with $\partial L\subset \partial B_{r_0}$,
$B_{\epsilon r_0}\cap L\ne \emptyset$
\begin{align}
\sup_{B_{r_0}\cap L} |A_L|^2 &\leq
\epsilon^2 \, r_0^{-2} \, ,\label{e:oo2}\\
\sup_{B_{r_0}\cap \Sigma} |A|^2 &\leq \epsilon_0^2 \, r_0^{-2}\, ,\label{e:oo3}
\end{align}
and $B_{r_0}\cap \Sigma$ lies on one side of $L$, then
$|A(z)|^2 \leq C \, \epsilon^2 \, r_0^{-2}$.
\end{Lem}

\begin{proof}
We will prove this for
$B_{r_0} = B_{r_0}(0) \subset \RR^n$
($z=0$, $\rc = \infty$); the general
case is similar (cf. \cite{CM1}).
Choose $\epsilon_0 > 0$ so:
If
$\cB_{2s}(y) \subset \Sigma$, $s \, \sup_{\cB_{2s}(y)} |A| \leq 4 \,
\epsilon_0$, and $t \leq 9s/5$, then $\in \Sigma_{y,t}$
is a graph over $T_y \Sigma$ with gradient $\leq t/s$ and
\begin{equation} \label{e:gmap1}
\inf_{y' \in \cB_{2s}(y)} |y'-y| / \dist_{\Sigma}(y,y') > 9/10 \, .
\end{equation}
Using $B_{\epsilon r_0}\cap L\ne \emptyset$,
let $L_{\frac{r_0}{2}}$ be a component of $B_{\frac{r_0}{2}} \cap L$
containing some $y_L \in B_{\epsilon r_0} \cap L$.
By \eqr{e:oo2} and \eqr{e:gmap1},
$L_{\frac{r_0}{2}} \subset \cB_{\frac{3r_0}{4}}(y_L)$.
Hence, by \eqr{e:oo2},
we can  rotate so $L_{\frac{r_0}{2}}$  is a graph over $\{ x_n =0 \}$
with $|\nabla_{L} x_n| \leq  \epsilon$ and
$ |x_n (L) | \leq 4 \, \epsilon \, r_0$.
Since $L \cap \Sigma = \emptyset$,   $x_n + 4
\, \epsilon \, r_0 >0$ is harmonic on $\cB_{\frac{r_0}{4}}
\subset \Sigma$.
By \eqr{e:oo3}, the Harnack inequality (and $0 \in \Sigma$) yields
$C=C(n)$ so
\begin{equation}        \label{e:hiy}
        0 < \sup_{\cB_{ \frac{r_0}{6} }} ( x_n + 4 \epsilon \, r_0 ) \leq
C  \, \inf_{\cB_{ \frac{r_0}{6} }} (x_n + 4
\, \epsilon \, r_0)
\leq 4 \, C \, \epsilon \,
r_0 \, .
\end{equation}
Since  $\cB_{\frac{r_0}{2}}$ is a graph with bounded gradient,
elliptic estimates give
\begin{equation} \label{e:oo6}
\int_{ \cB_{\frac{r_0}{8} }} |A|^2
\leq C' \, r_0^{-4} \, \int_{ \cB_{\frac{r_0}{6} }}
|x_n|^2 \, ,
\end{equation}
where $C' = C'(n)$. Combining \eqr{e:hiy} and
\eqr{e:oo6},   the lemma follows from the mean value inequality
since $\Delta |A|^2\geq -2\,|A|^4$ (see \cite{CM1}).
\end{proof}

\begin{Thm} \label{t:t2.1}
(Curvature estimate). There exist $\ec (i_0,k,n) ,
\rc (i_0,k,n)>0$ so: If $r_0 \leq \rc$,
$\Sigma^{n-1}\subset B_{r_0}=B_{r_0}(x)\subset M^n$ is an
embedded minimal hyper-surface, $\partial \Sigma \subset
\partial B_{r_0}$, and $\Sigma$ satisfies the
$(\ec,r_0)$-one-sided Reifenberg condition at $x$, then
$\sup_{B_{r_0-\sigma}\cap \Sigma}|A|^2\leq \sigma^{-2}$ for
$0< \sigma \leq r_0$.
\end{Thm}

\begin{proof}
Take $\rc > 0$ as in Lemma \ref{l:l2.1}, and set
$F=(r_0-r)^2\, |A|^2$.
Since $F \geq 0$, $F|\partial B_{r_0}\cap \Sigma=0$,  and $\Sigma$ is
compact, $F$ achieves its supremum at $y \in \partial
B_{r_0-\sigma}\cap \Sigma$ with $0< \sigma \leq r_0$. If $F\leq
1$, the theorem follows trivially. Hence, we may suppose
$F(y)=\sup_{B_{r_0}\cap \Sigma} F> 1$.
With $\epsilon_0\leq 1$ as in Lemma \ref{l:l2.1}, define $s>0$ by
$s^2\,|A(y)|^2=\epsilon_0^2/4$.
Since $F(y)=\sigma^2 \, |A(y)|^2> 1$ and $\epsilon_0 \leq 1$, we have
$2s< \sigma$. Since $F(y) >1$,
\begin{equation} \label{e:o2.3}
\sup_{B_s(y)\cap \Sigma}
\left( \frac{\sigma}{2} \right) ^2\, |A|^2
\leq \sup_{B_{\frac{\sigma}{2}}(y)\cap \Sigma}
\left( \frac{\sigma}{2} \right) ^2\, |A|^2
\leq \sup_{B_{\frac{\sigma}{2}}(y)\cap \Sigma} F
= \sigma^2\, |A(y)|^2\, .
\end{equation}
Multiplying \eqr{e:o2.3} by $4\,s^2/\sigma^2$ gives
$\sup_{B_{s}(y)\cap \Sigma} s^2\, |A|^2\leq \epsilon_0^2$.
Hence, the $(\ec , r_0)$-one-sided Reifenberg
assumption, Lemma \ref{l:l2.1} contradicting the choice of $s$ if $C \,
\ec^2 < \epsilon_0^2 /4$. Therefore, $F\leq 1$ for this $\ec$, and the
theorem follows.
\end{proof}

Letting $r_0 \to \infty$
in Theorem \ref{t:t2.1} gives the Bernstein type result:

\begin{Cor}
There exists $\epsilon (n)>0$ so any connected properly
embedded minimal hyper-surface satisfying the
$(\epsilon,\infty)$-one-sided Reifenberg condition is a hyper-plane.
\end{Cor}

We close by giving a condition which implies the
one-sided Reifenberg condition. Its proof (left to the reader)
relies on a simple
barrier argument (as in the proof of Corollary \ref{c:barrier}).

\begin{Lem} \label{l:l2.3}
There exist $\epsilon_0 (i_0,k) , r_1
(i_0,k)>0$, $c(i_0,k)\geq 1$ so: Let $\Sigma^{2}\subset
B_{r_0}=B_{r_0}(x)\subset M^{3}$ be an embedded minimal
disk, $\partial \Sigma \subset \partial B_{r_0}$, and
$r_0\leq r_1$. If for some $\epsilon<\epsilon_0$,
all $\sigma < r_0$ and all  $y \in
B_{r_0-\sigma} \cap \Sigma$ there is a minimal surface
$\Sigma_{y,\sigma}\subset B_{\sigma}(y) \setminus \Sigma$
with $\partial \Sigma_{y,\sigma} \subset \partial B_{\sigma}(y)$
and  $\Sigma_{y, \sigma} \cap B_{\epsilon\,\sigma}(y)\ne
\emptyset$, then $\Sigma$ satisfies the $(c \,
\epsilon_,r_0)$-one-sided Reifenberg condition at $x$.
\end{Lem}

\section{Laminations}   \label{s:lamination}

A codimension one {\it lamination}
on a $3$-manifold $M^3$ is a collection $\cL$ of
smooth disjoint surfaces (called leaves) such that
$\cup_{\Lambda \in \cL} \Lambda$ is closed.
Moreover, for each $x\in M$ there exists an
open neighborhood $U$ of $x$ and a coordinate chart, $(U,\Phi)$, with
$\Phi (U)\subset \RR^3$
so in these coordinates the leaves in $\cL$
pass through $\Phi (U)$ in slices of the
  form $(\RR\times \{ t\})\cap \Phi(U)$.
A {\it foliation} is a lamination for which
the union of the leaves is all of $M$
and a {\it minimal lamination} is a lamination whose leaves are minimal.
Finally, a sequence of laminations is said to converge
if the corresponding coordinate maps converge.
Note that any (compact) embedded surface (connected or not) is a lamination.

\begin{Pro}     \label{p:lamination}
Let $M^3$ be a fixed $3$-manifold.  If
$\cL_i\subset B_{2R}(x)\subset M$ is a sequence of minimal
laminations
with uniformly bounded curvatures (where each leaf has boundary contained
in $\partial B_{2R}(x)$), then a subsequence, $\cL_j$,
converges in the $C^{\alpha}$ topology for any $\alpha<1$ to
a (Lipschitz) lamination $\cL$ in $B_{R}(x)$ with minimal
leaves.
\end{Pro}

\begin{proof}
For convenience we will assume that each lamination $\cL_i$ has only
finitely many leaves where the number of leaves may depend on $i$.
This is all that is needed in the application of this proposition anyway.
Fix $x_0\in B_R(x)$.  The proposition will follow
once we construct uniform coordinate charts $\Phi_i$ on a
ball $B_{r_0}=B_{r_0}(x_0)$, where $4r_0\leq R$ is to be chosen.

By assumption, there exists $C$ so that
$\sup_{B_{4r_0} \cap \Lambda} |A|^2\leq C\, r_0^{-2}$ for each $i$
and every $\Lambda \in \cL_i$.
Replacing $r_0>0$ with a smaller radius,
 we may assume that
$C>0$ and $r_0\, \sqrt{k}$ are
as small as we wish and
$r_0< \frac{i_0}{2}$ ($i_0$ being the injectivity radius and $k$ a
bound for the curvature of $M$ in $B_{4r_0}$).
  In fact, if
$(x_1 , x_2 , x_3)$ are exponential normal coordinates centered at
$x_0$ on $B_{r_0}$, then
$\cup_{\Lambda\in \cL_i} B_{r_0}\cap \Lambda$ gives a sequence of
disconnected small curvature surfaces in these coordinates.
By standard estimates for normal coordinates,
the curvature is also small
with respect to the Euclidean metric.
Going to a further
subsequence (possibly with $r_0$ even smaller),
for each $i$ every sheet of
$\cup_{\Lambda\in \cL_i} B_{2r_0}(0)\cap \Lambda$ is a
graph with small gradient over a subset of the $\RR^2\times \{0\}$ plane
containing a ball of radius $r_0$ centered at the origin.

We claim that, in this ball, the sequence of laminations
converges in the $C^{\alpha}$ topology to a lamination for any $\alpha<1$.
The coordinate chart $\Phi$ required by the definition of a lamination will be
given by the
Arzela-Ascoli theorem as a limit of
a sequence of
bi-Lipschitz maps $\Phi_i:(x_j)_{j}\to \RR^3$ with bi-Lipschitz constants
close to one
and
defined on a slightly smaller concentric ball $B_{sr_0}$ for some
$s>0$ to be determined.
 Furthermore, we will show that for each $i$ fixed
$\Phi_i^{-1}(B_{sr_0}\cap \cup_{\Lambda\in\cL_i}\Lambda)$
is the union of planes which
are each parallel to $\RR^2\times \{0\}\subset \RR^3$; cf. \cite{So}.

Set  $\Phi_i^{-1}(y_1,y_2,y_3)=(y_1,y_2,\phi_i(y_1,y_2,y_3))$,
where $\phi_i$ is
defined as follows:
Order the sheets of
$B_{2r_0}(0)\cap_{\Lambda\in \cL_i} \Lambda$
as $\Lambda_{i,k}$ for $k=1,\cdots$ by
increasing values of $x_3$ and let $\Lambda_{i,k}$ be the graph of the
function $f_{i,k}$ over (part of) the $\RR^2\times \{ 0\}$ plane.
The domain of $f_{i,k}$ contains the ball of radius $r_0$ around the origin
in the $\RR^2\times \{ 0\}$ plane.  With slight abuse of notation,
 we will also denote balls
in $\RR^2\times \{ 0\}$ with
radius $t$ and center $0$ by $B_{t}$.
Set $w_{i,k} = f_{i,k+1}-f_{i,k}$.
In the next equation, $\Delta$, $\nabla$,
and $\text{div}$ will be  with respect to
the Euclidean metric on $\RR^2\times \{0\}$.
By a standard computation (cf.
\cite{Si},  (7) on p. 333
or chapter 1 of \cite{CM1}),
\begin{equation}    \label{e:endeq}
    \Delta w_{i,k}
    =  \text{div} \, (a\, \nabla w_{i,k} )
    +b \, \nabla w_{i,k}  + c \, w_{i,k}  \, .
\end{equation}
Here $a$ is a matrix valued function, $b$ is a vector valued function and
$c$ is a function.  Although $a$, $b$, and $c$ depend on $i$,
their scale invariant norms
are small when $C$ and $\sqrt{k} \, r_0$ are.
By \eqr{e:endeq}, the Schauder estimates and Harnack inequality
(e.g., 6.2 and 8.20
of \cite{GiTr})
applied to the positive function
$w_{i,k}$  give
\begin{equation}  \label{e:endhar}
    s r_0 \, \sup_{B_{s \, r_0}}\, |\nabla w_{i,k} |
    \leq
    C \, \sup_{B_{2s \, r_0}}\, w_{i,k}
    \leq
        \exp (\ccg\, s^{\beta} )\, \inf_{B_{2s \, r_0}}\, w_{i,k} \, .
\end{equation}
Where $\ccg$ and $\beta > 0$ depend on
the scale invariant norms of
$a,b$, and $c$.
Set  $\M_{i,k}=f_{i,k} (0)$.
In the region
$\{(y_1,y_2,y_3)\in B_{r_0}\times [\M_{i,k},\M_{i,k+1}]\}$,
define $\phi_i$ by
\begin{equation}
    \phi_i(y_1,y_2,y_3)
            =f_{i,k}(y_1,y_2)
          + \frac{y_3-\M_{i,k}}{\M_{i,k+1}-\M_{i,k}}
           w_{i,k}(y_1,y_2) \, .
\end{equation}
Hence
\begin{equation}   \label{e:o8.1end2}
    \nabla \phi_i = \nabla f_{i,k}
           +\frac{y_3-\M_{i,k}}{\M_{i,k+1}-\M_{i,k}}
           \nabla w_{i,k}
           +\frac{w_{i,k}}{\M_{i,k+1}-\M_{i,k}}
           \,\frac{\partial}{\partial y_3}\, .
\end{equation}
It follows easily from \eqr{e:endhar} and \eqr{e:o8.1end2} that for each $i$
the map $\Phi_i$ restricted to $B_{sr_0}(0)\subset \RR^3$ is bi-Lipschitz with
bi-Lipschitz constant close to one if $s$ is sufficiently small.
By the Arzela-Ascoli theorem, a
subsequence of $\Phi_i$ converges
in the $C^{\alpha}$ topology for any $\alpha < 1$
to a Lipschitz coordinate chart
$\Phi$ with the properties that are required.
The leaves in $B_{r_0}$ are $C^{1, \alpha}$
limits of minimal graphs with
bounded gradient, and hence minimal by elliptic regularity.
\end{proof}

Trivial examples show that the Lipschitz regularity above is optimal.

\section{A standard consequence of the maximum principle} \label{s:p4}

Using the maximum principle
and the convexity of small extrinsic balls we can bound the
topology of the intersection of a minimal surface with a ball:

\begin{Lem} \label{l:l1.2}
Let
$\Sigma^2\subset M^n$ be
an immersed minimal surface, $\partial \Sigma \subset
\partial B_{r_0}(x)$, $\K_{M^n}\leq k$, and injectivity radius of $M$
$\geq i_0$.  If
$r_0 < \min \{ \frac{i_0}{4} , \frac{\pi}{4 \, \sqrt{k} } \} $,
$B_t (y) \subset B_{r_0} (x)$,
and $\gamma \subset B_{t}(y) \cap \Sigma $
is a closed one-cycle homologous to zero in $B_{r_0} (x) \cap \Sigma$,
then $\gamma$ is homologous to zero in $B_{t}(y) \cap \Sigma$.
\end{Lem}

\begin{proof}
Apply the maximum principle to $f = \dist_M^2 \, (y,\cdot)$ on
the $2$-current that $\gamma$ bounds.
\end{proof}

By Lemma \ref{l:l1.2},
if $y\in B_t(x) \cap \Sigma$
is connected, then
$\chi (B_s(y)\cap \Sigma)
\geq \chi (B_t(x)\cap \Sigma)$ for $s+\dist_M(x,y)\leq t <
\min \{ \frac{i_0}{4} , \frac{\pi}{4 \, \sqrt{k} } \} $.
(The Euler characteristic is monotone.)

\section{A generalization of Proposition \ref{p:twosheets}}

In \cite{CM6}, the next proposition is needed
 when we deal with
the analog of the genus one helicoid (cf. \cite{HoKrWe}) where
$\Sigma$ (as above \eqr{e:defE1}) is not a disk. The {\it genus} of a surface
$\Sigma$ ($\Genus (\Sigma)$) is the genus of the closed surface
given by adding a disk to each boundary circle.

\begin{Pro} \label{p:twosheetsg}
There exist $C_0 ,  \epsilon_0$ so if  $0 \in \Sigma$, $\partial
\Sigma$ is connected,
 $\Genus (\Sigma) = \Genus (\Sigma_{0,r_1})$,
$R \geq C_0 \, r_1$, and $\epsilon_0 \geq \epsilon$, then $E_1
\cap \Sigma \setminus \Sigma_1$ is an (oppositely oriented)
$N$-valued graph $\Sigma_2$.
\end{Pro}

\begin{proof}
Note that, by the maximum principle and elementary topology (as in
part I of \cite{CM5}), $\Sigma \setminus \Sigma_{0,t}$ is an
annulus for $r_1 \leq t < 4R$. The proof now follows that of
Proposition \ref{p:twosheets}.

  First, (a slight extension of)
the ``estimate between the sheets'' given in theorem III.2.4 of
\cite{CM3} gives $\epsilon_0$ so that $E_1 \cap \Sigma$ is locally
graphical (this extension uses that $\Sigma \setminus
\Sigma_{0,r_1}$ is an annulus instead of that $\Sigma$ is a disk;
the proof of this extension is outlined in appendix A of
\cite{CM8}). As before, we get the second (oppositely oriented)
multi-valued graph $\Sigma_2 \subset \Sigma$.

Second, we argue by contradiction to show that there are no other
components of $E_1 \cap \Sigma$. Fix $\sigma_1 , \sigma_2$ as
before. The proof of Lemma \ref{l:onecomp} applies virtually
without change (since at least one of $\Sigma_a , \Sigma_b$ must
be a disk), so $\sigma_1$ and $\sigma_2$ connect in $\Sigma_{0,C_s
r_1}$. Hence,  $\sigma_0 \subset \partial \Sigma_{0,C_s r_1}$
connects $\sigma_1$ and $\sigma_2$.  Replace $\sigma_i$ with
$\sigma_i \setminus B_{C_s r_1}$, so that $\sigma_0 \cup \sigma_1
\cup \sigma_2 \subset \Sigma \setminus \Sigma_{0,r_1}$ is a simple
curve and
 $\partial
(\sigma_0 \cup \sigma_1 \cup \sigma_2) \subset
 \partial \Sigma_{0,R}$.
Let $\hat{\Sigma}$ be the component of $\Sigma_{0,R} \setminus
(\sigma_0 \cup \sigma_1 \cup \sigma_2)$ which does not intersect
 $\Sigma_{0,r_1}$.  It follows that $\hat{\Sigma}$ has genus zero and
connected boundary; i.e., it is a disk. Solve as above for the
stable disk $\Gamma$ with $\partial \Gamma = \partial
\hat{\Sigma}$ so $\Gamma$
 contains two disjoint $(N/2 -1)$-valued
graphs in $E_1$ which spiral together.
  For $R/r_1$ large,
Proposition \ref{c:0.4ofcm5} gives the point of large curvature,
contradicting the curvature estimate for stable surfaces.
\end{proof}

\end{document}